\documentclass[english]{article}
\usepackage[T1]{fontenc}
\usepackage[utf8]{inputenc}
\usepackage{geometry}
\usepackage{babel}
\usepackage{verbatim}
\usepackage{mathtools}
\usepackage{bm}
\usepackage[ruled,vlined,linesnumbered]{algorithm2e}
\usepackage{amsmath}
\usepackage{amsthm}
\usepackage{amssymb}
\usepackage{color}
\usepackage{graphicx}
\usepackage[numbers]{natbib}
\usepackage{multirow}
\usepackage{varwidth}
\providecommand{\tabularnewline}{\\}

\usepackage[pdfusetitle,
bookmarks=true,bookmarksnumbered=false,bookmarksopen=false,
 breaklinks=false,pdfborder={0 0 1},backref=false,colorlinks=false]
 {hyperref}

\global\long\def\lb{{\rm lb}}%

\makeatletter
\theoremstyle{plain}
\newtheorem{theorem}{Theorem}[section]
\theoremstyle{plain}
\newtheorem{proposition}[theorem]{\protect\propositionname}
\theoremstyle{plain}
\newtheorem{corollary}[theorem]{\protect\corollaryname}
\theoremstyle{definition}

\theoremstyle{remark}
\newtheorem{remark}{Remark}[section] 

\usepackage{babel}
\usepackage{bm}
\usepackage{etoolbox}
\usepackage{booktabs}

\renewcommand{\frac}[2]{\tfrac{#1}{#2}}

\let\oldsum\sum
\renewcommand{\sum}{\mathop{\textstyle\oldsum}\limits}

\newcommand{\apm}{\textrm{APMM}}

\global\long\def\APM{\textrm{APMM}}
\global\long\def\PMM{\textrm{PMM}}
\newcommand{\ifix}{\textrm{Inexact Fixed Point Method}}
\newcommand{\isecant}{\textrm{Inexact Secant Method}}

\providecommand{\corollaryname}{Corollary}
\providecommand{\propositionname}{Proposition}
\providecommand{\definitionname}{Definition}

\makeatother

\providecommand{\corollaryname}{Corollary}
\providecommand{\definitionname}{Definition}
\providecommand{\propositionname}{Proposition}

\begin{document}
\global\long\def\P{{\rm P}}%
\global\long\def\D{\textrm{{\rm D}}}%
\global\long\def\inprod#1#2{\left\langle #1,#2\right\rangle }%
\global\long\def\inner#1#2{\langle#1,#2\rangle}%
\global\long\def\binner#1#2{\big\langle#1,#2\big\rangle}%
\global\long\def\Binner#1#2{\Big\langle#1,#2\Big\rangle}%

\global\long\def\norm#1{\lVert#1\rVert}%
\global\long\def\abs#1{\left |#1\right |}%
\global\long\def\bnorm#1{\big\Vert#1\big\Vert}%
\global\long\def\Bnorm#1{\Big\Vert#1\Big\Vert}%

\global\long\def\setnorm#1{\Vert#1\Vert_{-}}%
\global\long\def\bsetnorm#1{\big\Vert#1\big\Vert_{-}}%
\global\long\def\Bsetnorm#1{\Big\Vert#1\Big\Vert_{-}}%

\global\long\def\brbra#1{\big(#1\big)}%
\global\long\def\Brbra#1{\Big(#1\Big)}%
\global\long\def\rbra#1{(#1)}%
\global\long\def\sbra#1{[#1]}%
\global\long\def\bsbra#1{\big[#1\big]}%
\global\long\def\Bsbra#1{\Big[#1\Big]}%
\global\long\def\cbra#1{\{#1\}}%
\global\long\def\bcbra#1{\big\{#1\big\}}%
\global\long\def\Bcbra#1{\Big\{#1\Big\}}%

\global\long\def\vertiii#1{\left\vert \kern-0.25ex  \left\vert \kern-0.25ex  \left\vert #1\right\vert \kern-0.25ex  \right\vert \kern-0.25ex  \right\vert }%

\global\long\def\matr#1{\bm{#1}}%

\global\long\def\til#1{\tilde{#1}}%
\global\long\def\wtil#1{\widetilde{#1}}%

\global\long\def\wh#1{\widehat{#1}}%

\global\long\def\mcal#1{\mathcal{#1}}%

\global\long\def\mbb#1{\mathbb{#1}}%

\global\long\def\mtt#1{\mathtt{#1}}%

\global\long\def\ttt#1{\texttt{#1}}%

\global\long\def\dtxt{\textrm{d}}%

\global\long\def\bignorm#1{\bigl\Vert#1\bigr\Vert}%
\global\long\def\Bignorm#1{\Bigl\Vert#1\Bigr\Vert}%

\global\long\def\rmn#1#2{\mathbb{R}^{#1\times#2}}%

\global\long\def\deri#1#2{\frac{d#1}{d#2}}%
\global\long\def\pderi#1#2{\frac{\partial#1}{\partial#2}}%

\global\long\def\limk{\lim_{k\rightarrow\infty}}%

\global\long\def\smid{\mskip1mu\mid\mskip1mu}%

\global\long\def\trans{\textrm{T}}%

\global\long\def\onebf{\mathbf{1}}%
\global\long\def\zerobf{\mathbf{0}}%
\global\long\def\qbf{\mathbf{q}}%
\global\long\def\zero{\mathbf{0}}%


\global\long\def\Euc{\mathrm{E}}%
\global\long\def\Expe{\mathbb{E}}%

\global\long\def\rank{\mathrm{rank}}%
\global\long\def\range{\mathrm{range}}%
\global\long\def\diam{\mathrm{diam}}%
\global\long\def\epi{\mathrm{epi} }%
\global\long\def\relint{\mathrm{relint} }%
\global\long\def\dom{\mathrm{dom}}%
\global\long\def\prox{\mathrm{prox}}%
\global\long\def\proj{\mathrm{Proj}}%
\global\long\def\for{\mathrm{for}}%
\global\long\def\diag{\mathrm{diag}}%
\global\long\def\and{\mathrm{and}}%
\global\long\def\var{\textrm{VaR}}%
\global\long\def\where{\mathrm{where}}%
\global\long\def\dist{\mathrm{dist}}%

\global\long\def\True{\mathbf{True}}%
\global\long\def\False{\mathbf{False}}%

\global\long\def\st{\mathrm{s.t.}}%

\global\long\def\inte{\operatornamewithlimits{int}}%
\global\long\def\cov{\mathrm{Cov}}%

\global\long\def\argmin{\operatornamewithlimits{arg\,min}}%
\global\long\def\argmax{\operatornamewithlimits{arg\,max}}%
\global\long\def\maximize{\operatornamewithlimits{maximize}}%
\global\long\def\minimize{\operatornamewithlimits{minimize}}%

\global\long\def\tr{\operatornamewithlimits{tr}}%

\global\long\def\dis{\operatornamewithlimits{dist}}%

\global\long\def\prob{{\rm Pr}}%

\global\long\def\spans{\textrm{span}}%
\global\long\def\st{\operatornamewithlimits{s.t.}}%
\global\long\def\subjectto{\textrm{subject to}}%

\global\long\def\Var{\operatornamewithlimits{Var}}%

\global\long\def\raw{\rightarrow}%
\global\long\def\law{\leftarrow}%
\global\long\def\Raw{\Rightarrow}%
\global\long\def\Law{\Leftarrow}%

\global\long\def\vep{\varepsilon}%

\global\long\def\dom{\operatornamewithlimits{dom}}%

\global\long\def\tsum{{\textstyle {\sum}}}%

\global\long\def\Cbb{\mathbb{C}}%
\global\long\def\Ebb{\mathbb{E}}%
\global\long\def\Fbb{\mathbb{F}}%
\global\long\def\Nbb{\mathbb{N}}%
\global\long\def\Rbb{\mathbb{R}}%
\global\long\def\Sbb{\mathbb{S}}%
\global\long\def\reals{\mathbb{R}}%

\global\long\def\extR{\widebar{\mathbb{R}}}%
\global\long\def\Pbb{\mathbb{P}}%

\global\long\def\Mrm{\mathrm{M}}%
\global\long\def\Acal{\mathcal{A}}%
\global\long\def\Bcal{\mathcal{B}}%
\global\long\def\Ccal{\mathcal{C}}%
\global\long\def\Dcal{\mathcal{D}}%
\global\long\def\Ecal{\mathcal{E}}%
\global\long\def\Fcal{\mathcal{F}}%
\global\long\def\Gcal{\mathcal{G}}%
\global\long\def\Hcal{\mathcal{H}}%
\global\long\def\Ical{\mathcal{I}}%
\global\long\def\Kcal{\mathcal{K}}%
\global\long\def\Lcal{\mathcal{L}}%
\global\long\def\Mcal{\mathcal{M}}%
\global\long\def\Ncal{\mathcal{N}}%
\global\long\def\Ocal{\mathcal{O}}%
\global\long\def\Pcal{\mathcal{P}}%
\global\long\def\Scal{\mathcal{S}}%
\global\long\def\Tcal{\mathcal{T}}%
\global\long\def\Wcal{\mathcal{W}}%
\global\long\def\Xcal{\mathcal{X}}%
\global\long\def\Ycal{\mathcal{Y}}%
\global\long\def\Zcal{\mathcal{Z}}%


\global\long\def\abf{\mathbf{a}}%
\global\long\def\bbf{\mathbf{b}}%
\global\long\def\cbf{\mathbf{c}}%
\global\long\def\dbf{\mathbf{d}}%
\global\long\def\ebf{\mathbf{e}}%
\global\long\def\fbf{\mathbf{f}}%
\global\long\def\gbf{\mathbf{g}}%
\global\long\def\pbf{\mathbf{p}}%
\global\long\def\sbf{\mathbf{s}}%
\global\long\def\rbf{\mathbf{r}}%
\global\long\def\lbf{\mathbf{l}}%
\global\long\def\ubf{\mathbf{u}}%
\global\long\def\vbf{\mathbf{v}}%
\global\long\def\wbf{\mathbf{w}}%
\global\long\def\xbf{\mathbf{x}}%
\global\long\def\ybf{\mathbf{y}}%
\global\long\def\zbf{\mathbf{z}}%

\global\long\def\Abf{\mathbf{A}}%
\global\long\def\Cbf{\mathbf{C}}%
\global\long\def\Ubf{\mathbf{U}}%
\global\long\def\Pbf{\mathbf{P}}%
\global\long\def\Ibf{\mathbf{I}}%
\global\long\def\Ebf{\mathbf{E}}%
\global\long\def\Mbf{\mathbf{M}}%
\global\long\def\Qbf{\mathbf{Q}}%
\global\long\def\Lbf{\mathbf{L}}%
\global\long\def\Pbf{\mathbf{P}}%
\global\long\def\Wbf{\mathbf{W}}%

\global\long\def\lambf{\bm{\lambda}}%
\global\long\def\mubf{\bm{\mu}}%
\global\long\def\alphabf{\bm{\alpha}}%
\global\long\def\sigmabf{\bm{\sigma}}%
\global\long\def\thetabf{\bm{\theta}}%
\global\long\def\deltabf{\bm{\delta}}%
\global\long\def\vepbf{\bm{\vep}}%
\global\long\def\pibf{\bm{\pi}}%


\global\long\def\abm{\bm{a}}%
\global\long\def\bbm{\bm{b}}%
\global\long\def\cbm{\bm{c}}%
\global\long\def\dbm{\bm{d}}%
\global\long\def\ebm{\bm{e}}%
\global\long\def\fbm{\bm{f}}%
\global\long\def\gbm{\bm{g}}%
\global\long\def\hbm{\bm{h}}%
\global\long\def\pbm{\bm{p}}%
\global\long\def\qbm{\bm{q}}%
\global\long\def\rbm{\bm{r}}%
\global\long\def\sbm{\bm{s}}%
\global\long\def\tbm{\bm{t}}%
\global\long\def\ubm{\bm{u}}%
\global\long\def\vbm{\bm{v}}%
\global\long\def\wbm{\bm{w}}%
\global\long\def\xbm{\bm{x}}%
\global\long\def\ybm{\bm{y}}%
\global\long\def\zbm{\bm{z}}%

\global\long\def\Abm{\bm{A}}%
\global\long\def\Bbm{\bm{B}}%
\global\long\def\Cbm{\bm{C}}%
\global\long\def\Dbm{\bm{D}}%
\global\long\def\Ebm{\bm{E}}%
\global\long\def\Fbm{\bm{F}}%
\global\long\def\Gbm{\bm{G}}%
\global\long\def\Hbm{\bm{H}}%
\global\long\def\Ibm{\bm{I}}%
\global\long\def\Jbm{\bm{J}}%
\global\long\def\Lbm{\bm{L}}%
\global\long\def\Obm{\bm{O}}%
\global\long\def\Pbm{\bm{P}}%
\global\long\def\Qbm{\bm{Q}}%
\global\long\def\Rbm{\bm{R}}%
\global\long\def\Sbm{\bm{S}}%
\global\long\def\Ubm{\bm{U}}%
\global\long\def\Vbm{\bm{V}}%
\global\long\def\Wbm{\bm{W}}%
\global\long\def\Xbm{\bm{X}}%
\global\long\def\Ybm{\bm{Y}}%
\global\long\def\Zbm{\bm{Z}}%
\global\long\def\lambm{\bm{\lambda}}%

\global\long\def\alphabm{\bm{\alpha}}%
\global\long\def\albm{\bm{\alpha}}%
\global\long\def\pibm{\bm{\pi}}%
\global\long\def\taubm{\bm{\tau}}%
\global\long\def\mubm{\bm{\mu}}%
\global\long\def\yrm{\mathrm{y}}%
\global\long\def\ifapl{\texttt{ifapl}}

\newcommand{\zhenwei}[1]{{#1}}
\newcommand{\blue}[1]{{#1}}

\newcommand{\red}[1]{\textcolor{red}{#1}}

\global\long\def\aleq{\overset{(a)}{\leq}}%

\global\long\def\bleq{\overset{(b)}{\leq}}%

\global\long\def\cleq{\overset{(c)}{\leq}}%

\global\long\def\dleq{\overset{(d)}{\leq}}%

\global\long\def\ageq{\overset{(a)}{\geq}}%

\global\long\def\bgeq{\overset{(b)}{\geq}}%

\global\long\def\cgeq{\overset{(c)}{\geq}}%

\global\long\def\beq{\overset{(b)}{=}}%

\global\long\def\ceq{\overset{(c)}{=}}%

\global\long\def\deq{\overset{(d)}{=}}%

\global\long\def\vbfp{\vbf_{\text{p}}}%

\global\long\def\vbfd{\vbf_{\text{d}}}%

\newcommand{\remove}[1]{\textcolor{gray}{[\textbf{remove:} #1]}} 

\global\long\def\lb{{\rm lb}}%

\title{Uniformly Optimal and Parameter-free First-order Methods for Convex and Function-constrained
Optimization}

\author{
Qi Deng\thanks{Antai College of Economics \& Management, Shanghai Jiao Tong University. Email: qdeng24@sjtu.edu.cn} 
$\qquad$  Guanghui Lan \thanks{
Industrial and Systems Engineering, Georgia Institute of Technology. Email: george.lan@isye.gatech.edu
}
$\qquad$  Zhenwei Lin\thanks{School of Industrial Engineering, Purdue University. Email: lin2193@purdue.edu} 
}

\maketitle
\begin{abstract}

This paper presents new first-order methods for achieving optimal oracle complexities in convex optimization with convex function constraints. Oracle complexities are measured by the number of function and gradient evaluations. To achieve this, we develop first-order methods that can utilize computational oracles for solving diagonal quadratic programs in subproblems.
For problems where the optimal value  $f^*$  is known, such as those in overparameterized models and feasibility problems, we propose an accelerated first-order method that incorporates a modified Polyak step size and Nesterov's momentum. Notably, our method does not require knowledge of smoothness levels, H\"{o}lder continuity parameter of the gradient, or additional line search, yet achieves the optimal oracle complexity bound of  $\mathcal{O}(\varepsilon^{-2/(1+3\rho)})$ under H\"{o}lder smoothness conditions.
When  $f^*$  is unknown, we reformulate the problem as finding the root of the optimal value function and develop inexact fixed-point iteration and secant method to compute  $f^*$. These root-finding subproblems are solved inexactly using first-order methods to a specified relative accuracy. We employ the accelerated prox-level (APL) method, which is proven to be uniformly optimal for convex optimization with simple constraints. Our analysis demonstrates that APL-based root-finding also achieves the optimal oracle complexity of $\mathcal{O}(\varepsilon^{-2/(1+3\rho)})$ for convex function-constrained optimization, without requiring knowledge of any problem-specific structures. Through experiments on various tasks, we demonstrate the advantages of our methods over existing approaches in function-constrained optimization. 

\end{abstract}

\section{Introduction}
In this paper, we are interested in solving the following nonlinear programming problem:
\begin{equation}
{f^*}=\min_{\xbf \in\Xcal} \ f(\xbf),\ \ 
\st\ \ g_{i}(\xbf)\le0,\quad i\in [m],
\label{pb:func-constraint}
\end{equation}
where $f(\xbf)$ and $g_{i}(\xbf)$ are convex continuous functions, and $\Xcal\subseteq\reals^{d}$ is a closed convex set, nonempty and typically polyhedral. Notation $[m]$ is short for $\{1,2,\ldots, m\}$. 
Both \( f(\xbf) \) and \( g_i(\xbf) \) can be either nonsmooth, smooth or weakly smooth, with the level of smoothness unknown a priori. Convex optimization with inequality constraints, as formulated in \eqref{pb:func-constraint}, has seen a resurgence of interest in fields such as machine learning and operations research.
This interest stems from numerous applications, such as Neyman-Pearson classification~\citep{rigollet2011neyman}, risk-averse learning~\citep{cheng2022functional}, and fairness in machine learning~\citep{cotter2019optimization}, among others. 
Traditionally, interior point methods have been known for solving problem~\eqref{pb:func-constraint} with high accuracy. However, these methods are inefficient for large-scale problems due to the repetitive need to solve Newton's system. To address large-scale problems, first-order methods~\citep{hamedani2021primal, xu2021iteration, boob2019proximal}, which bypass the need to compute the Hessians of $f$ and $g$, have become the primary tool and have attracted significant research attention.

Among these works, one popular approach is the penalty method, including the augmented Lagrangian method (e.g. \citep{LanMon13-1,xu2021iteration}), which repeatedly applies first-order methods to inexactly solve the regularized problem penalized by the constraint violation. 
For instance, \citet{xu2021iteration} proposed an inexact augmented Lagrangian method that employs Nesterov's accelerated gradient to solve the proximal subproblem. This method achieves an $\mathcal{O}(1/\varepsilon)$ complexity bound for smooth convex-constrained optimization and improves the bound to ${\Ocal}(1/\sqrt{\varepsilon}\log(1/\vep))$ when the objective function is strongly convex. 
\citet{hamedani2021primal} extended the renowned primal-dual hybrid gradient method \citep{chambolle2011first} to more general convex-concave saddle point problems, including convex function-constrained optimization~\eqref{pb:func-constraint} as a special case. To address the challenge of unbounded domains, they incorporated a backtracking line search into the primal-dual method, achieving a complexity bound of $\Ocal(1/{\varepsilon})$ for solving problem~\eqref{pb:func-constraint}. \citet{lin2022efficient} developed a new accelerated primal-dual method to deal with strongly convex constraint functions. By progressively leveraging the strong convexity of the Lagrangian function, they obtained an improved complexity bound of ${\Ocal}(1/\sqrt{\varepsilon})$.
 A unified constraint extrapolation method (ConEx, \citet{boob2019proximal}) has been developed to handle both stochastic and deterministic problems, as well as convex and strongly convex problems. 
Another important direction for solving constrained optimization is the level set method~\cite{aravkin2019level, lin2018level, nesterov2018lectures}, which reformulates the original problem as finding the root of a convex function. Specifically, \citet{lin2018level} proposed a feasible level-set algorithm that applies the fixed point iteration to search the optimal value $f^*$. In this method, the max-type subproblem can be approximately solved using smooth approximation, followed by applying Nesterov's accelerated method.

Despite recent advances, the achieved \emph{oracle complexities}, which are measured in terms of \emph{function value and gradient evaluations} in this paper, remain inferior to the lower bounds of first-order methods in convex optimization. Particularly, it is well-established that the lower complexity bound for smooth convex optimization is $\mathcal{O}(1/\sqrt{\varepsilon})$, which can be attained by Nesterov's accelerated gradient method~\citep{nesterov1983method}. However, standard studies in first-order methods typically assume the existence of an abstract convex set domain onto which projection operators have closed-form solutions. This assumption becomes impractical in scenarios involving nonlinear constraints, highlighting the greater challenges inherent in function-constrained optimization. 

To further reduce oracle complexity, it seems reasonable to sacrifice computational efficiency by employing stronger projection oracles.  
\citet[Sec 2.3.5]{nesterov2018lectures} described a new constrained minimization scheme that can achieve the optimal oracle complexity. Unlike standard first-order methods, this approach requires solving a structured quadratically constrained quadratic program (QCQP) where the quadratic term has a diagonal structure:
\begin{equation}\label{eq:qcqp}
    \min_{\xbf\in\Xcal} \quad \norm{\xbf-\bar\xbf}^2 \ \  
    \st \quad \norm{\xbf-\zbf_i}^2 \le b_i, \ i\in[m].
\end{equation}
Later, \citet{boob2024level} introduced a feasible level-constrained first-order method that similarly solves a diagonal QCQP. They obtained the optimal oracle complexity to reach first-order stationary points in smooth nonconvex optimization and showed that the diagonal QCQP can often be efficiently solved using either first-order methods or customized interior point methods.
\citet{zhang2022solving} proposed an accelerated constrained gradient descent (ACGD) method, which achieves the optimal oracle complexity bounds in smooth and strongly convex optimization while only solving a relatively easier diagonal quadratic program (QP):
\begin{equation}\label{eq:lcqp}
    \min_{\xbf\in\Xcal} \quad \norm{\xbf-\bar\xbf}^2 \ \ 
    \st \quad \abf_i^\top \xbf \le b_i, \ i\in[m].
\end{equation}
As pointed out by  \citet[Sec. 2.4]{devanathan2023polyak}, such a special quadratic problem can often be solved efficiently by utilizing the diagonal structure. 
Despite these recent advances in convex function-constrained optimization, several limitations remain. First, all these methods are developed for Lipschitz smooth optimization and cannot be readily applied when the smoothness level is unknown. Second, they involve multiple parameters that require careful fine-tuning. In particular, they require knowing the smooth parameter $L$ to determine the stepsizes. Even when the parameter $L$ is known, additional line search is often necessary, as the dual domain and the Lipschitz parameter of the Lagrangian function are unbounded~\citep{zhang2022solving,hamedani2021primal}.

The search for efficient parameter-free methods has been a long-standing research area, which, to the best of our knowledge, largely involves no functional constraints. By “parameter-free,” we mean that all required parameters can be set independently of specific problem instances, without influencing the convergence or complexity bounds.
\citet{polyak1987introduction} introduced an adaptive stepsize which has been found to be more efficient than simple diminishing stepsizes in the subgradient method~\citep{boyd2003subgradient}. Polyak's stepsize is particularly useful for solving nonlinear equations or over-parameterized models~\citep{loizou2021stochastic}, where the optimal value  $f^*$  is known (often zero or near zero in these cases). \citet{hazan2019revisiting} observed that Polyak's stepsize can automatically adapt to the problem structure. However, their result does not achieve optimal convergence rates in the convex smooth setting.
Inspired by Polyak's work, \citet{devanathan2023polyak} considered convex function-constrained optimization where $f^*$ is known.  A related approach is the level bundle method~\citep{kiwiel1995proximal,LNN, Ben-Tal05, de2014level}, which has primarily focused on nonsmooth problems. 
\citet{van2014level, tang2022new} developed restricted memory level bundle methods for convex and nonsmooth function-constrained optimization. 
\citet{lan2015bundle} extended the bundle-level method to smooth optimization and proposed accelerated variants that are uniformly optimal under different levels of H\"{o}lder smoothness. 
For a comprehensive review of bundle-type methods, we refer to \citet{frangioni2020standard}.
While the bundle-level method~\citep{lan2015bundle} requires solving a series of quadratic problems, \citet{nesterov2015universal} introduced line search-based gradient methods that are universally optimal under different smoothness levels and solve easier subproblems.  Line search-free and parameter-free gradient methods that achieve optimal rates have been investigated in several recent works~\citep{li2023simple, rodomanovuniversal, lan2024projected}.

\subsection{Contributions}
We develop new algorithms for function-constrained convex optimization that achieve optimal oracle complexity bounds for various smoothness levels, requiring limited or no parameter tuning.

First, we introduce an accelerated Polyak minorant method~(\apm{}) for solving convex and function-constrained optimization when the optimal value $f^*$ is known. This algorithm achieves the optimal oracle complexity $\mathcal{O}(\left({1}/{\varepsilon}\right)^{{2}/{(1+3\rho)}})$ under H\"{o}lder smoothness conditions, where $\rho$ denotes the smoothness level. Remarkably, it requires no prior knowledge of $\rho$ or other smoothness parameters, except for the optimal value $f^*$. The algorithm solves a structured quadratic problem with a diagonal quadratic term~\eqref{eq:lcqp}. Furthermore, we develop a restarted variant of \apm{} that achieves faster convergence rates under certain H\"{o}lder error bound conditions.
\apm{} is closely related to Polyak's stepsize method and, more broadly, to bundle-level methods~\citep{polyak1987introduction, devanathan2023polyak, Ben-Tal05, lan2015bundle} for convex optimization. Specifically, it is an \emph{accelerated bundle method with a fixed bundle level}. \blue{Unlike the accelerated prox-level (APL) method~\citep{lan2015bundle}, which requires a fixed proximal center and employs a double-loop procedure to search for the optimal level, \apm{} uses a moving proximal center and operates as a single-loop algorithm.} When applied to unconstrained convex problems, \apm{} can be viewed as a variant of Polyak-type stepsize method accelerated by Nesterov's momentum. 

Second, when the optimal value $f^*$ is unknown, we propose infeasible level-set methods to identify $f^*$. Specifically, we formulate the function-constrained problem~\eqref{pb:func-constraint} as an equivalent convex root-finding problem associated with the optimal value function~\citep{devanathan2023polyak,lin2018level} 
$V(\eta):\Rbb\rightarrow\Rbb$, which is defined as: \[
V(\eta)\coloneqq \min_{\xbf\in\Xcal}v(\xbf,\eta),\text{where }v(\xbf,\eta)\coloneqq \max\bcbra{f(\xbf)-\eta,g_{1}(\xbf),\ldots,g_{m}(\xbf)}.
\]
Note that when $\eta=f^*$, we have $V(f^*)=\min_{\xbf\in\Xcal}v(\xbf,f^*)=0$.
\blue{We solve the root-finding problem using an inexact variant of the fixed point iteration, which achieves linear convergence to the optimal value $f^*$. To further enhance efficiency, we introduce a novel truncated secant stepsize strategy that combines the advantages of both fixed and secant stepsizes. This approach enables the use of the more aggressive secant stepsize when appropriate, while reverting to the more stable fixed stepsize in the presence of significant inexactness. Our empirical results demonstrate the practical effectiveness of this strategy.}
The level-set subproblem, which is nonsmooth due to the max-type structure of $v(\eta,\xbf)$, can be approximately solved by the APL method~\citep{lan2015bundle}. 
\blue{We introduce techniques to further improve practical performance. Since APL is invoked repeatedly within the root-finding procedure, having a good initial gap estimate can greatly accelerate convergence. We address this by developing a warm-start strategy that provides sharp gap estimates, thereby speeding up the APL method.}
We show that the total oracle complexity of APL-based fixed point iteration retains the $\mathcal{O}(\norm{\ybf^*}_1+1)\log(\norm{\ybf^*}_1+1)\left({1}/{\varepsilon}\right)^{{2}/{(1+3\rho)}})$ bound where $\ybf^*$ is the vector of Lagrange multipliers. This result is further improved to $\Ocal(\min\{(\norm{\ybf^*}_1+1)\log(\norm{\ybf^*}_1+1),\log (1/\vep)\}\left({1}/{\varepsilon}\right)^{{2}/{(1+3\rho)}})$ for the APL-based secant method, which is insensitive to the poor conditioning for large $\norm{\ybf^*}_1$. Notably, both approaches do not require prior knowledge of $\norm{\ybf^*}_1$ or parameter-tuning. 
To the best of our knowledge, this paper establishes the first optimal oracle complexity results for general convex function-constrained optimization under H\"{o}lder smoothness conditions. A detailed comparison of the complexity results is provided in Table~\ref{tab:complexity-compare}.

Third, we conduct experiments on a variety of constrained problems, including solving the KKT systems of Second-Order Cone Programming and Linear Matrix Inequalities, convex Quadratically Constrained Quadratic Programming, and Neyman-Pearson classification problems. Numerical results are quite encouraging, showing that our methods exhibit strong performance on both nonsmooth and smooth problems.
While bundle methods have traditionally been considered most effective for nonsmooth optimization, our study suggests that they also hold significant potential for solving smooth and function-constrained optimization problems.

\begin{table}
\small
\caption{\label{tab:complexity-compare}Oracle complexities for convex function-constrained optimization,  which are measured by the number of function and gradient evaluations. The parameter $L$ refers to the known Lipschitz constant of the gradient.  PJ denotes the standard projection on $\Xcal$. }
\centering{}%
\begin{tabular}{ccccc}
\hline 
\multirow{2}{*}{Methods} & Function & Known & Subproblem & \multirow{2}{*}{ Complexities}\tabularnewline
 & types & parameters & types & \tabularnewline
\hline 
LCPG~\citep{boob2024level} & Lipschitz smooth & $L$ & D-QCQP~\eqref{eq:qcqp} & $\mcal O(1/\vep)$\tabularnewline
\cline{1-1}
\multirow{1}{*}{ACGD~\citep{zhang2022solving}} & Lipschitz smooth & $L$ & D-QP~\eqref{eq:lcqp} & $\mcal O(1/\sqrt{\vep})$ \tabularnewline
\cline{1-1}
iALM~\citep{xu2021iteration} & Lipschitz smooth & $L$ & PJ & $\mcal O(1/\vep)$\tabularnewline
\cline{1-1}
\apm{} (this paper) & H\"{o}lder smooth & $f^{*}$ known & D-QP~\eqref{eq:lcqp} & $\mcal O((1/\vep)^{2/(1+3\rho)})$\tabularnewline
APL-based methods (this paper)
 & H\"{o}lder smooth & - & D-QP~\eqref{eq:lcqp} & ${\mcal O}((1/\vep)^{2/(1+3\rho)})$\tabularnewline
\hline 
\end{tabular}
\end{table}

\subsection{Preliminaries}\label{subsec:prelim}

Throughout this paper, we use bold letters for the vectors, such as $\xbf$ and $\ybf$. We use $\Rbb^d$ for the Euclidean space and $\Sbb^m$ for the  standard simplex: $\Sbb^m=\{\xbf\in\Rbb^{m+1}:\tsum_{i=0}^m x_i =1, x_i\ge0, 0\le i\le m\} $. 
$\norm{\cdot}$ and $\norm{\cdot}_{\infty}$ stand for the Euclidean norm and infinity norm, respectively. The indicator function and normal cone are given by $\iota_{\cal X}(\xbf)=0$ if $\xbf \in \mcal X$ and $\infty$ otherwise and $\mcal N_{\cal X}(\xbf)=\{\zbf:\inner{\zbf}{ \xbf-\ybf}\geq 0,\forall \ybf\in \mcal X\}$, respectively.
We assume that $f(\xbf)$ and $g_i(\xbf)$ are real-valued convex functions. Furthermore, we assume the existence of a real-valued mapping  $M_i:\reals^{d}\times\reals^{d}\rightarrow(0,\infty)$, ($0\le i \le m$) such that
\begin{equation}\label{eq:H_smooth}
\begin{aligned}
f(\ybf)-f(\xbf)-\inner{f'(\xbf)}{\ybf-\xbf} & \le\frac{M_{0}(\ybf,\xbf)}{1+\rho}\norm{\ybf-\xbf}^{1+\rho},\\
g_{i}(\ybf)-g_{i}(\xbf)-\inner{g'_{i}(\xbf)}{\ybf-\xbf} & \le\frac{M_{i}(\ybf,\xbf)}{1+\rho}\norm{\ybf-\xbf}^{1+\rho},
\end{aligned}    
\end{equation}
for any $\xbf,\ybf\in\Rbb^d$, where $\rho\in[0,1]$.  Depending on the value of  $\rho$, we categorize $f(\xbf)$ as follows: when  $\rho = 0$,  $f(\xbf)$  is nonsmooth; when  $\rho \in (0,1)$, $f(\xbf)$ is weakly smooth; and when  $\rho = 1$, $f(\xbf)$  is smooth.
For brevity, we use $f'(\xbf)$ to denote a subgradient of $f(\xbf)$ and the gradient when $f(\xbf)$ is differentiable. 
We frequently use the linear minorant, which is given by $\ell_{f}(\xbf,\ybf)\coloneqq f(\ybf)+\inner{f'(\ybf)}{\xbf-\ybf}$, and the maximum of minorants: $v_\ell(\xbf,\bar\xbf, \eta)\coloneqq \max\{\ell_f(\xbf, \bar\xbf)-\eta, \ell_{g_1}(\xbf,\bar{\xbf}),\ell_{g_2}(\xbf,\bar{\xbf}),\ldots, \ell_{g_m}(\xbf,\bar{\xbf})\}$.
We denote $\gbf(\xbf)=[g_{1}(\xbf),g_{2}(\xbf),\ldots,g_{m}(\xbf)]^{\top}\in\Rbb^{m}$ for brevity.  We say that a solution $\xbf\in\Xcal$ is $\vep$-optimal for problem~\eqref{pb:func-constraint} if $f(\xbf)-f^* \le \vep$ and $\norm{[\gbf(\xbf)]_+}_\infty \le \vep$, where $[\xbf]_+$ computes the positive parts element-wise.  

\subsection{Paper Structure}
Section~\ref{sec:apmm} introduces the accelerated Polyak’s minorant method for solving the constrained problem~\eqref{pb:func-constraint} when  $f^*$  is known. Section~\ref{sec:root-finding} discusses the general root-finding approach for solving \eqref{pb:func-constraint} when  $f^*$  is unknown. Section~\ref{sec:apl-based} presents more practical algorithms, where the accelerated prox-level method is applied to solving the root-finding subproblem. Finally, Section~\ref{sec:numerical} provides a variety of numerical applications to demonstrate the advantages of our proposed methods.

\section{The Accelerated Polyak's Minorant Method}\label{sec:apmm}
In this section, we consider solving problem~\eqref{pb:func-constraint} when $f^*$ is known. 
Our goal is to introduce a new accelerated gradient method that not only achieves the optimal oracle complexity but also eliminates the need for additional parameter tuning. The proposed method, referred to as \apm{}, is outlined in Algorithm~\ref{alg:APMM}.

Conceptually, \apm{} is similar to the classic accelerated gradient method (e.g., \citep[Alg. 1]{tseng2008accelerated}, \citep[Sec. 3.3]{lan2020first}), which involves three intertwined sequences, \( \{\xbf^k, \ybf^k, \zbf^k\} \). The primary distinction between the standard accelerated gradient method and \apm{} lies in the update rule for \( \{\xbf^k\} \).
To illustrate this difference, let's temporarily ignore the function constraint $ \gbf(\xbf) \le \zero$ in problem~\eqref{pb:func-constraint} and assume that $f(\xbf)$  is  $L$-Lipschitz smooth. In this scenario, the update of  $\xbf^k$ in accelerated methods typically has the form
\begin{equation}\label{eq:nag-update}
\xbf^k = \argmin_{\xbf \in \Xcal} \bsbra{ \inner{f'(\zbf^{k-1})}{\xbf - \zbf^{k-1}} + \frac{1}{\eta_k} \norm{\xbf - \xbf^{k-1}}^2}.
\end{equation}
The choice of stepsize $\eta_k$ usually depends on the Lipschitz parameter. A typical setting to ensure the $\Ocal(1/\sqrt{\varepsilon})$ complexity bound is $\eta_k={k}/{L}$ and $\alpha_k=2/(k+1)$. However, this approach requires knowledge of the curvature and the level of smoothness.
In contrast to the update in \eqref{eq:nag-update}, \apm{} introduces an additional cut constraint using the linear minorant while minimizing a stabilizing term associated with the previous iterate, $\xbf^{k-1}$. This allows \apm{} to bypass the need for knowing problem parameters other than $f^*$. Moreover,  \apm{} imposes a dynamically changing constraint $X_k$ in the update. While $X_k$ can be set to $\Xcal$ for simplicity, it can also be formulated using cutting planes over the past iterations to construct refined minorants. 
\begin{algorithm}[h]
\small
\KwIn{$\xbf^0$, $f^*$, $\vep$} 
Set $\ybf^{0}=\xbf^{0}$, $X_{0}=\Xcal$, $\bar{v}_0=v(\ybf^0,f^*)$\;
\For{$k=1,2,\ldots,$}{
Compute $\zbf^{k}, \xbf^k, \tilde{\ybf}^k$ such that:
\begin{equation}\label{eq:update-zk}
     \zbf^{k}=(1-\alpha_{k})\ybf^{k-1}+\alpha_{k}\xbf^{k-1},
\end{equation}
\begin{equation}\label{eq:update-xk}
    \xbf^{k}=\argmin_{\xbf\in X_{k-1}}  \, \frac{1}{2}\|\xbf-\xbf^{k-1}\|^2 \ 
\st\, v_\ell(\xbf,\zbf^k,f^*)\le 0 
\end{equation}
\begin{equation}
    \tilde{\ybf}^{k}=(1-\alpha_{k})\ybf^{k-1}+\alpha_{k}\xbf^{k}.
\end{equation}

\lIf{$\bar v_{k-1} < v(\tilde\ybf^k,f^*)$}{$\bar v_k = \bar v_{k-1},\ \ybf^k = \ybf^{k-1}$}
\lElse{$\bar v_k = v(\tilde\ybf^k,f^*),\ \ybf^k = \tilde\ybf^k$}
{\bf Break if } $\bar{v}_{k}\le\vep$\;
Update set $X_{k}$;} 
\KwOut{$\ybf^{k}$}
\caption{\textbf{A}ccelerated \textbf{P}olyak \textbf{M}inorant \textbf{M}ethod~(\apm{})}\label{alg:APMM}
\end{algorithm}


In the following theorem, we summarize the main convergence property of Algorithm~\ref{alg:APMM}.
\begin{theorem}\label{thm:rate-APM}
Let $\xbf^{*}$ be an optimal solution of problem~\eqref{pb:func-constraint} and $f^*=f(\xbf^*)$. Suppose $\xbf^*\in X_k$, $k=0,1,2,\ldots,$.
Let us define the sequence $\{\Gamma_{k}\}_{k\ge1}$
by $\Gamma_{1}=1$ and $\Gamma_{k}=\blue{\Pi_{i=2}^{k}(1-\alpha_{i})^{-1}}$ for $k\ge 2$, and $\cbf_{K}\coloneqq\bsbra{\alpha_{1}^{1+\rho}\Gamma_{1},\ldots,\alpha_{k}^{1+\rho}\Gamma_{k},\ldots,\alpha_{K}^{1+\rho}\Gamma_{K}}$ in Algorithm~\ref{alg:APMM}. After $K$ iterations, we have 
\[
v(\ybf^K, f^*) \le \blue{\Gamma_{K}^{-1}(1-\alpha_{1}) } v(\ybf^{0},f^{*})+\frac{\Gamma_{K}^{-1}\hat{M}}{1+\rho}\norm{\cbf_{K}}_{\frac{2}{1-\rho}}\|\xbf^{*}-\xbf^{0}\|^{\rho+1},
\]
where $\hat{M}_{}=\max_{0\le i\le m}\sup_{\bar{\xbf},\hat{\xbf}\in\{\xbf:\|\xbf^{*}-\xbf\|\le \|\xbf^{*}-\xbf^{0}\|\}}M_{i}(\bar{\xbf},\hat{\xbf})$. 
In particular, setting $\alpha_{k}=\frac{2}{k+1}$ gives
\begin{equation*}\begin{aligned}
\max \bcbra{f(\ybf^{K})-f^{*}, \norm{[\gbf(\ybf^{K})]_+}_\infty} & \le \frac{\hat{M}}{1+\rho}\|\xbf^{*}-\xbf^{0}\|^{\rho+1}\frac{2^{\rho+1}3^{(1-\rho)/2}}{K^{(1+3\rho)/2}}.
\end{aligned}\end{equation*}
\end{theorem}


\begin{remark}
    Our result indicates that Algorithm~\ref{alg:APMM} achieves the optimal oracle complexity~\citep{nesterov2015universal,lan2015bundle} across various smoothness levels. Notably, the only required parameter of Algorithm~\ref{alg:APMM} is the averaging scheme $\{\alpha_k\}$, which is problem-independent but crucial for achieving acceleration. For comparison, we will show that an alternative scheme corresponding to PMM~\citep{devanathan2023polyak} can establish a suboptimal complexity. Beyond setting $\alpha_k$, there is no need to know the smoothness level of the function gradient. Moreover, due to the stabilization property~\eqref{eq:sum-square-xk}, the H\"{o}lder smoothness condition only needs to be satisfied locally.
\end{remark}

\paragraph{Choice of $X_{k}$.} The choice of  $X_k$  is quite flexible, with the basic requirement being that it should include the optimal solutions to problem~\eqref{pb:func-constraint}.  We can set 
\begin{enumerate}
\item $X_{k}=\Xcal$. 
\item Full memory set: 
$X_{k}=X_{k-1}\cap\bcbra{\xbf:\ell_{f}(\xbf,\zbf^{k})\le f^{*}},\,k=1,2,\ldots,
X_{0} =\Xcal.$
\item A limited memory set: for some $k_0>0$, set
\[
X_{k}= \begin{cases}
 \Xcal \cap \bigcap_{1\le s\le k}\bcbra{\xbf:\ell_{f}(\xbf,\zbf^{s})\le f^{*}} & 1\le k<k_{0}\\
\Xcal \cap \bigcap_{k-k_{0}+1\le s\le k}\bcbra{\xbf:\ell_{f}(\xbf,\zbf^{s})\le f^{*}} & k\ge k_{0}
\end{cases}
\]
\item Averaging ${X}_{k} =\Xcal\cap\left\{ \xbf:\sum_{s=0}^{k}\beta_{s}\ell_{f}(\xbf,\zbf^{s})\le\left(\sum_{s=0}^{k}\beta_{s}\right)f^{*}\right\},$
where $\beta_{s}\ge0$, $s=0,1,\ldots,k$. 
\end{enumerate}

\paragraph{Connection with the Polyak's gradient method}
\blue{Suppose we drop the functional constraints (i.e., $m=0$).} Let $\Xcal=\Rbb^d$ and set $X_{k-1}=\Rbb^d$.  Then the subproblem in \eqref{eq:update-xk} reduces to
$\min_{\xbf\in\Rbb^{d}}\,\frac{1}{2}\norm{\xbf-\xbf^{k-1}}^{2}\;\st\ell_{f}(\xbf,\zbf^{k})\le f^{*}.$
Let $\Lcal(\xbf,\eta)=\eta[\ell_{f}(\xbf,\zbf^{k})-f^{*}]+\frac{1}{2}\norm{\xbf-\xbf^{k-1}}^{2}$
be the Lagrangian function. 
Applying the KKT condition yields
$\xbf^{k}=\xbf^{k-1}-\eta f'(\zbf^{k}), \ \text{where }\ \eta=\max\left\{ 0,\frac{\ell_{f}(\xbf^{k-1},\zbf^{k})-f^{*}}{\norm{f'(\zbf^{k})}^{2}}\right\}$.
If we set $\alpha_{k}=1$, then $\zbf^{k}=\xbf^{k-1}$, we have $\ell_{f}(\xbf^{k-1},\zbf^{k})-f^{*}=f(\xbf^{k-1})-f^{*}$. Consequently, \apm{} reduces to the Polyak's update $
\xbf^k = \xbf^{k-1} - \frac{f(\xbf^{k-1}) - f^*}{\norm{f'(\xbf^{k-1})}^2} f'(\xbf^{k-1}).$
Recently, \citet{devanathan2023polyak} proposed the Polyak Minorant Method (PMM), which is inspired by Polyak’s stepsize and can be extended to use various minorant functions. For simplicity, we focus on the case where the minorant is given by a linear approximation. In this context, their updating scheme reduces to a special case of our algorithm with  $\alpha_k = 1$.
\begin{theorem}\label{thm:pmm}
Let $\alpha_k=1$ in \apm{}. Then, it requires at most
\blue{
$K=\Big\lceil\left(\frac{\hat{M}}{1+\rho}\right)^{\frac{2}{1+\rho}}\norm{\xbf^0 - \xbf^*}^{\blue{2}}\cdot \vep^{-\frac{2}{1+\rho}}\Big\rceil$}
iterations to obtain an $\vep$-optimal solution.
\end{theorem}

\paragraph{Restart under H\"{o}lderian error bounds.}
We show that {\apm} exhibits even faster convergence rates under the  error-bound condition.  
Specifically, a continuous function $h:\Rbb^d\raw \Rbb$ has a \emph{H\"{o}lderian error bound} on $\Xcal$ if there exists $\mu,\tilde{\rho}>0$
such that for any $\xbf\in\Xcal$, the following holds: 
\begin{equation}\label{eq:H-growth}
h(\xbf)-h^{*}\geq\frac{\mu}{2}\norm{\xbf-\bar{\xbf}}^{\tilde{\rho}},\ \ \bar{\xbf}=\proj_{\Xcal^{*}}(\xbf),    
\end{equation}
where $h^*=\min_{\ybf\in\Xcal}h(\ybf)$ and $\Xcal^{*}:=\{\xbf\in\Xcal:h(\xbf)=h^*\}$. 
Condition~\eqref{eq:H-growth} is also known as the Hölderian growth condition, which coincides with the quadratic growth condition when $\tilde{\rho}=2$ (~\citep[Section 3.4]{Necoara2019Linear}) and reduces to sharpness when $\tilde{\rho}=1$~\citep{burke1993weak,davis2018subgradient}. The H\"{o}lderian growth condition has been extensively studied in various contexts~\cite{bolte2017error,jiang2022holderian}, such as the application to piecewise convex polynomial functions~\cite{li2013global}. We omit the case $1+\rho > \tilde{\rho}$ in our paper since it implies \blue{there exists one constant $c$ such that $\norm{\xbf-\xbf^*}\geq c>0,\forall \xbf\in \mcal X\backslash \bcbra{\xbf^*}$}, which does not seem intuitive. One example is when $\Xcal=\Xcal^*$.
Next, we develop a restarted scheme to achieve faster convergence under the growth condition.
\begin{algorithm}[h]
\small
\KwIn{$\blue{\bar{\xbf}^0}, f^*, \theta, \vep$;}
Set $s=0,\Delta_0=\max\{f(\blue{\bar{\xbf}^0})-f^{*},\norm{[\gbf(\blue{\bar{\xbf}^0})]_{+}}_{\infty}\}$\;
\While{$v(\blue{\bar{\xbf}^s},f^*)>\vep$}{
Compute $\blue{\bar{\xbf}^{s+1}}$=\apm{}($\blue{\bar{\xbf}^{s}},f^*, \Delta_0\cdot\theta^{s+1}$)\;
Set $s=s+1$;
}
\KwOut{$\blue{\bar{\xbf}^s}$}

\caption{Restarted {\apm} (r{\apm})}\label{alg:restartedAcc}
\end{algorithm}
\begin{theorem}\label{thm:restart_rate}
Suppose $f$ has a H\"{o}lderian growth on $\Xcal$ with $\mu,\tilde{\rho}>0$, $\theta\in(0,1)$. Let $\Delta_{0}=\max\{f(\blue{\bar{\xbf}^0})-f^{*},\norm{[\gbf(\blue{\bar{\xbf}^0})]_+}_{\infty}\}$, and  denote $\bcbra{\blue{\bar{\xbf}^{s}}}_{s\geq0}$
as the sequence generated by Algorithm~\ref{alg:restartedAcc}, then \blue{ for generating solution $\blue{\bar{\xbf}^{s+1}}$, \apm{} needs $K_{s+1}=\Big\lceil \mathfrak{C}\cdot\theta^{\frac{2(\rho+1-\tilde{\rho})s-2\tilde{\rho}}{\tilde{\rho}(1+3\rho)}}\Big\rceil$ oracle calls for generating $\bar{\xbf}^{s+1}$}
and the overall oracle complexity of Algorithm~\ref{alg:restartedAcc} to find an $\vep$-optimal solution is bounded by 
\begin{equation}
T_{\vep}=\begin{cases}
\mathfrak{C}\theta^{-\frac{2}{3\rho+1}}\cdot\Brbra{\theta^{\frac{2\rbra{\rho+1-\tilde{\rho}}}{\tilde{\rho}\rbra{1+3\rho}}}-1}^{-1}\cdot\blue{\brbra{\frac{\Delta_{0}}{\theta^{2}\vep}}^{\frac{-2\rbra{\rho+1-\tilde{\rho}}}{\tilde{\rho}\rbra{1+3\rho}}}}+2+\log_{1/\theta}\Brbra{\frac{\Delta_{0}}{\vep}}, & \text{if }\tilde{\rho} > 1+\rho, \\
\Brbra{\mathfrak{C}\theta^{-\frac{2}{3\rho+1}}+1}\cdot\Brbra{\log_{1/\theta}\Brbra{\frac{\Delta_{0}}{\vep}}+\blue{2}}, & \text{if }\tilde{\rho} = 1+\rho,
\end{cases}
\end{equation}
where $\mathfrak{C}:=\brbra{\frac{\hat{M}}{1+\rho}}^{\frac{2}{1+3\rho}}\cdot\brbra{\frac{2}{\mu}}^{\frac{2(\rho+1)}{\tilde{\rho}(1+3\rho)}}\cdot2^{\frac{2(\rho+1)}{1+3\rho}}\cdot3^{\frac{1-\rho}{1+3\rho}}\cdot \Delta_{0}^{\frac{2(\rho+1-\tilde{\rho})}{\tilde{\rho}(1+3\rho)}}$. 
\end{theorem}

\section{Root-finding for General Convex Function-constrained Problems}\label{sec:root-finding}

This section considers convex function-constrained optimization where $f^*$ is unknown. 
Searching for the optimal value ${f^*}$ can be cast as a root-finding problem: $
\text{Find} \ \blue{f^*} \coloneqq \min\{\eta: V(\eta)=0\}.$
Motivated by the study of level-set methods~\citep{lin2018level,aravkin2019level}, we propose to use root-finding algorithms, namely, the fixed-point iteration and the inexact secant method, to find the optimal level $\blue{f^*}$. 
The next section develops more specific algorithms where the subproblem is solved inexactly by the bundle-level methods.

Before presenting the main algorithms, we describe some useful properties of the value function. 
\begin{proposition}\label{prop:value-function}
The value function $V(\cdot)$ satisfies the following properties.
\begin{enumerate}
    \item $V:\reals \rightarrow \reals $ is convex, non-increasing. $\blue{f^*}=\min_\eta\{\eta:V(\eta)=0\}$.
    \item $V(\cdot)$ is  1-Lipschitz continuous, i.e., for $\Delta\ge0$,
$V(\eta)-\Delta\le V(\eta+\Delta)\le V(\eta).$
\item For any $\eta_1,\eta_2$ such that $\eta_1 < \eta_2$, we have 
$-1 \le V'(\eta_1)\le \frac{V(\eta_1)-V(\eta_2)}{\eta_1-\eta_2} \le V'(\eta_2) \le 0$.
\item Let $\bar{f}=\min_{\xbf\in\Xcal}f(\xbf)$.  We assume $\bar{f}<f^*$. Let $\bar\Xcal=\{\xbf\in\Xcal: f(\xbf)=\bar{f}\}$ and $\bar{g}=\min_{\xbf\in\bar\Xcal}\max_{1\le i \le m} g_i(\xbf)$, then $\bar{g}>0$.  In other words, all the solutions in $\bar\Xcal$ are infeasible for \eqref{pb:func-constraint}.  Moreover, for any $\eta \le \bar{f}-\bar{g}$, we have $V(\eta)=\bar{f}-\eta$. 
\item Suppose $\eta<\blue{f^*}$ and the assumption of part 4 holds, then any minimizer $\tilde\xbf\in \argmin_{\xbf\in\Xcal} v(\xbf,\eta)$ is an infeasible point of \eqref{pb:func-constraint}. 
\end{enumerate} 
\end{proposition}

\begin{remark}
The monotonicity, convexity, and Lipschitz continuity of $V(\cdot)$ are known from prior work~\citep{nesterov2018lectures, lin2018level}. The assumption in Part 4 of Proposition~\ref{prop:value-function} is mild, and it implies that the constraints are non-negligible. Even if this assumption is not satisfied, detecting such a degenerate case can be done efficiently. For further details, see our discussion in Section~\ref{subsec:initialization}.
\end{remark}

\subsection{The inexact fixed point iteration}

The above discussion motivates us to develop an infeasible level-set method for finding the optimal level $\blue{f^*}$. Specifically, we start from a sufficiently small initial value $\eta_0$ ($\eta_0<\blue{f^*}$) and apply a root-finding algorithm to generate a sequence $\{\eta_t\}$ converging to $\blue{f^*}$ \blue{from below}. Under the mild non-degeneracy condition~(Proposition~\ref{prop:value-function}), the optimal solutions $\xbf^t\in\argmin_\xbf v(\xbf,\eta_t)$, \blue{where $\eta_t < f^*$}, are infeasible for problem~\eqref{pb:func-constraint}.
However, classic root-finding algorithms require an exact evaluation of $V(\cdot)$, which is often impossible. To bypass this issue, we assume that the lower and upper bounds $l, u$ of $V(\eta)$ are available, such that 
\begin{equation}\label{eq:subprob-relative-error}
\frac{u}{l}\le\alpha,\quad 0< l\le V(\eta)\le u,
\end{equation}
for some constant $\alpha > 1$.
Next, we develop the inexact fixed point iteration in Algorithm~\ref{alg:Inexact-fixed-point} and establish its convergence properties in Theorem~\ref{thm:convergence-fixedpoint-abstract}.

\begin{algorithm}\label{alg:Inexact-fixed-point}
\small
\KwIn{$\eta_0<\blue{f^*}$, $\alpha>1$, $\beta\in(0,1)$, $\vep\in(0,\infty)$; 
}
Find $l_0, u_0$ such that $0< l_0\le V(\eta_0) \le u_0\le \alpha l_0$\;
\For{$t=1,2,\ldots,$}{
Compute $\eta_{t}=\eta_{t-1}+\beta l_{t-1}$\;
Find a couple $(l_{t},u_{t})$ such that \blue{$\frac{u_t}{l_t}\leq \alpha$, $0<l_t\leq V(\eta_t)\leq u_t$ and \blue{$u_t\leq u_{t-1}$}\;}
{\bf Break if }$u_{t}\le\vep$\;
}
\KwRet{$\eta$.}
\caption{Inexact Fixed Point method (with fixed stepsize)}
\end{algorithm}

\begin{theorem}\label{thm:convergence-fixedpoint-abstract}
Under the assumptions of Algorithm~\ref{alg:Inexact-fixed-point},
we have that $\eta_{t}$ is monotonically increasing and $\eta_{t}\le\blue{f^*}$
for all $t\ge0$. 
Moreover, \blue{let $\bar{V}^{\prime}=\min\bcbra{\xi:\xi\in \partial V(f^*)}$}, then we have 
\begin{equation}\label{eq:etak-recursive-1}
\blue{f^*}-\eta_{t}\le\sigma(\blue{f^*}-\eta_{t-1}),\quad t=1,2,\ldots, \text{ where } \sigma\coloneqq 1+\frac{\beta}{\alpha} \blue{\bar{V}^{\prime}}.
\end{equation}
Then the algorithm terminates in at most 
\blue{$T^{\rm{FP}}_\vep=\Big\lceil\frac{\alpha}{-\beta {\bar{V}^{\prime}}}\log\left(\frac{\alpha V(\eta_0)}{-{\bar{V}^{\prime}}\vep}\right)\Big\rceil$}
iterations. 
\end{theorem}

\paragraph{The condition number}
The above analysis implies that the convergence rate of fixed point iteration depends on the condition number $1/|\blue{\bar{V}^{\prime}}|$. {A similar inexact fixed point iteration was proposed by \cite{lin2018level}, which is initiated from the right side of the root: $\eta_0>\blue{f^*}$. Their condition number is related to the measure $(\eta_0-\blue{f^*})/|V(\eta_0)|$ , which appears to be worse than ours.} Next, we further exploit the connection between the subdifferential $\partial V(\blue{f^*})$ and the optimality condition. 
\begin{theorem}\label{thm:V-deri}
\blue{Suppose the KKT conditions for problem~\eqref{pb:func-constraint} are satisfied at an optimal solution (which, for instance, is ensured by Slater's condition). Then, we have}
$ 
    \partial V(\blue{f^*})\supseteq\big\{-\frac{1}{1+\norm{\ybf^*}_1}: \  \ybf^*\in\Rbb^m_+ \ \text{is a vector of Lagrange multipliers}\big\}.$
\end{theorem}

\begin{corollary}\label{cor:thm4_5_cor}
    \blue{Assume that the conditions of Theorems~\ref{thm:convergence-fixedpoint-abstract} and~\ref{thm:V-deri} are satisfied. Then, Algorithm~\ref{alg:Inexact-fixed-point} will terminate in at most $T^{\rm{FP}}_\vep=\Big\lceil\frac{\alpha\rbra{1+\norm{\ybf^{*}}_{1}}}{\beta}\log\left(\frac{\alpha V(\eta_{0})\rbra{1+\norm{\ybf^{*}}_{1}}}{\vep}\right)\Big\rceil$ iterations.}
\end{corollary}
\blue{Corollary~\ref{cor:thm4_5_cor} follows directly by combining Theorems~\ref{thm:convergence-fixedpoint-abstract} and~\ref{thm:V-deri}.}
According to Corollary~\ref{cor:thm4_5_cor}, the inexact fixed point iteration achieves a complexity of $\Ocal((1+\norm{\ybf^*}_1)\log(1/\vep))$, where the conditioning is determined by the Lagrange multipliers. 
In the following, we introduce a variant of the secant method whose performance is less sensitive to the conditioning.

\subsection{The inexact secant method}
The secant method approximates the derivative via finite differences: $\eta_t=\eta_{t-1}-\frac{\eta_{t-2}-\eta_{t-1}}{V(\eta_{t-2}) - V(\eta_{t-1})} V(\eta_{t-1})$,
which exhibits a superlinear rate of convergence to the optimum.
Intuitively, since the ratio $\frac{V(\eta_{t-2}) - V(\eta_{t-1})}{\eta_{t-2} - \eta_{t-1}}$ falls within the interval $[-1, 0)$, the secant method is more aggressive than a fixed stepsize.
As the optimal value $V(\eta)$ cannot be computed exactly,  \citet{aravkin2019level} considered the inexact secant method involving the finite-difference approximation: $
\eta_t=\eta_{t-1}-\frac{\eta_{t-2}-\eta_{t-1}}{u_{t-2}-l_{t-1}} l_{t-1}.$

 However, due to the inexactness in computing $V(\eta)$, the resulting stepsize can be more conservative when $\{u_t, l_t\}$ are inaccurate. Hence, we develop a new truncated secant stepsize by integrating the fixed-point iteration and the traditional secant method, as outlined in Algorithm~\ref{alg:inexact-secant}. We establish the convergence of  Algorithm~\ref{alg:inexact-secant} in the following theorem, and the proof of this result is similar to that of \cite{aravkin2019level}. 
 \begin{algorithm}[h]
 \small 
 \caption{Truncated Inexact Secant method}\label{alg:inexact-secant}
\KwIn{$\alpha$, $\beta$, $\vep$, $\eta_0, \eta_1$ such that $\eta_0<\eta_1<\blue{f^*}$;}

Obtain $(l_{0},u_{0},\eta_{0})$
and $(l_{1},u_{1},\eta_{1})$\;

\For{$t=2,3,\ldots,$}{
Compute $\eta_{t}=\eta_{t-1}+\beta\max\{1,-\frac{\eta_{t-2}-\eta_{t-1}}{u_{t-2}-l_{t-1}}\}l_{t-1}$\; \label{tag:eta-update-secant}
Find a couple $(l_{t},u_{t})$ such that \blue{$\frac{u_t}{l_t}\leq \alpha,$ $0<l_t\leq V(\eta_t)\leq u_t$ and $u_t\leq u_{t-1}$\;}
{\bf Break if }$u_{t}\le\vep$\;
}
\end{algorithm}

\begin{theorem} \label{thm:convergence-secant-1} 
Suppose that we set $\beta\in(1/2,1]$ and $\alpha\in (1,2\sqrt{\beta})$. Then, $\{\eta_t\}$ is a monotonically increasing sequence with $\eta_t\le \blue{f^*}$ for all $t>0$. Moreover,  the contraction property~\eqref{eq:etak-recursive-1} holds. \blue{Let $\bar{V}^{\prime}=\min\bcbra{\xi:\xi\in\partial V\brbra{f^{*}}}$}. Suppose we obtain $\eta_1$ from $\eta_0$ by one step of fixed point iteration, then the total iteration number of Algorithm~\ref{alg:inexact-secant} for $\vep$-optimal solution is bounded by
\blue{
$T^{\rm{SC}}_\vep=\Big\lceil\min\Big\{ {{\frac{\alpha}{-\beta  \blue{\bar{V}^{\prime}}}\log\brbra{\frac{\alpha V(\eta_0)}{-\blue{\bar{V}^{\prime}}\vep}}},\log_{2\sqrt{\beta}/\alpha}\brbra{\frac{\alpha l_0}{\varepsilon} \sqrt{\frac{\blue{f^*}-\eta_0}{\eta_1-\eta_0}}}}\Big\}\Big\rceil.$}
\end{theorem}

\begin{remark}
Theorem~\ref{thm:convergence-secant-1} implies that Algorithm~\ref{alg:inexact-secant} exhibits two convergence patterns. The linear rate driven by the fixed stepsize has an appealing contraction property, which means each $\eta_t$ is getting closer to $\blue{f^*}$ by a fixed ratio.  However, this rate is influenced by the conditioning of the problem, determined by $1/|\blue{\bar{V}^{\prime}}|$.
On the other hand, due to the secant-type stepsize, Algorithm~\ref{alg:inexact-secant} also exhibits a more robust linear rate, which is independent of $1/|\blue{\bar{V}^{\prime}}|$.
\end{remark}

\blue{We further give the complexity bound with respect to the dual variable in the following corollary.}
\begin{corollary}
    \blue{Assume that the conditions of Theorem~\ref{thm:V-deri} and~\ref{thm:convergence-secant-1} are satisfied. Then, Algorithm~\ref{alg:inexact-secant} will terminate in $T^{\rm{SC}}_\vep=\Big\lceil\min\Big\{{{\frac{\alpha\brbra{1+\norm{\ybf^{*}}_{1}}}{\beta}\log\brbra{\frac{\alpha V(\eta_{0})\brbra{1+\norm{\ybf^{*}}_{1}}}{\vep}}},\log_{2\sqrt{\beta}/\alpha}\brbra{\frac{\alpha l_{0}}{\varepsilon}\sqrt{\frac{\blue{f^{*}}-\eta_{0}}{\eta_{1}-\eta_{0}}}}}\Big\}\Big\rceil$ iterations.}
\end{corollary}

\section{Root-finding based on Bundle-level Methods}\label{sec:apl-based}
In this section, we present more concrete root-finding procedures, ensuring that the subproblem is solved efficiently to meet the condition defined in \eqref{eq:subprob-relative-error}. 
A critical component is to evaluate $V(\eta)$, which can be framed as the following problem:
\begin{equation}\label{pb:value-func}
\min_{\xbf\in\Xcal}\:v(\xbf,\eta)\coloneqq\max\bcbra{f(\xbf)-\eta,g_{1}(\xbf),\ldots,g_{m}(\xbf)}.
\end{equation}
We apply the accelerated prox-level method (APL)~\citep[Sec. 3]{lan2015bundle}, which is both parameter-free and capable of generating verifiable optimality gaps~\eqref{eq:subprob-relative-error}. Moreover, APL obtains the optimal rates on the H\"{o}lder smooth and convex problems. The convergence of APL requires a bounded domain assumption. Therefore, throughout the rest of the paper, we assume that the domain is bounded, i.e., $D_\Xcal \coloneqq \max_{\xbf,\ybf\in \mcal X}\norm{\xbf-\ybf}<\infty$.
The rest proceeds as follows. 

\begin{algorithm}[h]\label{alg:gap-reduce}
\small
\KwIn{$\text{\ensuremath{\xbf}},\tilde{l},\eta$, $\theta\in(0,1)$;}
$\ybf^0=\xbf^{0}=\xbf$, $\tilde{u}=v_0^{{\rm U}}=v(\xbf^{0},\eta)$,
$v_{0}^{{\rm L}}=\tilde{l}$, $\bar{\Xcal}_{0}=\Xcal$, $\lambda=\frac{1}{2}(v_0^{{\rm L}}+v_0^{{\rm U}})$, $k=1$\;

\While{True}{
Update $\zbf^{k}$ by $\zbf^{k}=(1-\alpha_k)\ybf^{k-1}+\alpha_k \xbf^{k-1}$\;
Compute $h_{k} =\min_{\xbf\in\bar{\Xcal}_{k-1}} v_{\ell}(\xbf,\zbf^k,\eta)$ and $v_{k}^{{\rm L}} =\max\{v_{k-1}^{{\rm L}},\min\{\lambda,h_{k}\}\}$\;
\lIf{\blue{$\tilde{u}-v_{k}^{L}\le \frac{1+\theta}{2}(\tilde{u}-\tilde{l})$}}{\label{tag:gap-reduce-break-1} Set $\pbf=\ybf^{k-1}$, $\tilde{v}^{{\rm L}}=v_{k}^{{\rm L}}$, and \KwSty{break}}

Compute $\xbf^{k}$ by solving
\begin{equation}\label{eq:update-xk-2}
    \min_{\xbf\in\bar{\Xcal}_{k-1}} \ \|\xbf-\xbf^0\|^2 \quad \quad 
\st  \ v_{\ell}(\xbf,\zbf^{k},\eta)\le \lambda.
\end{equation}

Choose a set $\bar{\Xcal}_{k}$ such that $\Xcal_{k}^{{\rm L}}\subseteq\bar{\Xcal}_{k}\subseteq\Xcal_{k}^{{\rm U}}$, where $\Xcal_{k}^{{\rm L}}=\bcbra{\xbf\in\bar{\Xcal}_{k-1}:\, v_{\ell}(\xbf,\zbf^{k},\eta)\le \lambda},
\Xcal_{k}^{{\rm U}}=\bcbra{\xbf\in\Xcal:\inner{\xbf^k-\xbf^0}{\xbf-\xbf^k}\ge 0 }.$

Update $\tilde{\ybf}^{k}=\alpha_{k}\xbf^{k}+(1-\alpha_{k})\ybf^{k-1}$\;
Compute $v_{k}\coloneqq v(\tilde{\ybf}^{k},\eta)$; 


\lIf{${v}^{\rm U}_{k-1} < v_{k}$}{${v}_{k}^{\rm U}= {v}^{\rm U}_{k-1}\ ,\ \ybf^k = {\ybf}^{k-1}$}
\lElse{${v}_{k}^{\rm U}= v_{k}\ ,\ \ybf^k = \tilde{\ybf}^{k}$}

\lIf{$\blue{v_{k}^{U}-\tilde{l}\le \frac{1+\theta}{2}(\tilde{u}-\tilde{l})}$}{\label{tag:gap-reduce-break-2} Set $\pbf=\ybf^{k}$, $\tilde{v}^{{\rm L}}=v_{k}^{{\rm L}}$, and \KwSty{break}}
Set $k\leftarrow k+1$\;
}
\KwOut{$(\pbf, \tilde{v}^{\rm L}$);}

\caption{Gap reduction $\Gcal(\xbf,\tilde{l},\eta, \theta)$}

\end{algorithm}

\subsection{APL for solving the subproblem in root-finding}\label{subsec:APL}
The APL method is motivated by the classic bundle-level method~\citep{LNN, Ben-Tal05} and further uses Nesterov's momentum to obtain faster convergence rates. 
A key component is a subroutine called gap reduction, which is adapted for solving \eqref{pb:value-func} and described in Algorithm~\ref{alg:gap-reduce}. In essence, Algorithm~\ref{alg:gap-reduce} maintains the lower bound $v_k^{\rm L}$ and upper bound $v_k^{\rm U}$ of $V(\eta)$, iteratively refining $v_k^{\rm L}$ through linear minorants and $v_k^{\rm U}$ using an accelerated scheme similar to \apm{}. Since the target value $V(\eta)$ is unknown, Algorithm~\ref{alg:gap-reduce} employs a guessed value $\lambda$. It is possible that the optimal solution is not feasible within the set constraints, and therefore, the algorithm explicitly requires a bundle constraint \( \bar{\Xcal}_{k-1} \). However, the extra constraint can involve as few as one more cut constraint provided we choose $\Xcal_k^{\rm U}$, which does not substantially increase the computation burden. 
Algorithm~\ref{alg:gap-reduce} reduces the optimality gap by a constant ratio: $\tilde{v}_k^{\rm{U}} - \tilde{v}_k^{\rm{L}} \le (\frac{1+\theta}{2})(\tilde{u}-\tilde{l})$, which is determined by the input parameter $\theta$. 
Convergence analyses of Algorithm~\ref{alg:gap-reduce} are summarized in Appendix~\ref{sec:convergence_ana_sec_4}.
For more illustrations, we refer to \citet[Section 3]{lan2015bundle}.

\begin{algorithm}[h]\label{alg:apl}
\small
\caption{The \textbf{A}ccelerated \textbf{P}roximal \textbf{L}evel (APL) method, $\Acal(\xbf^0, \bar{l}_0, \eta, \theta, \alpha)$}
    \KwIn{($\xbf^0$, $\bar{l}_0$, $\eta$, $\alpha$);}
    Set $\bar{u}_0=v(\xbf^0,\eta)$, $s=0$\;
    \While{$\bar{u}_s-\bar{l}_s > \frac{\alpha-1}{\alpha} \bar{u}_s$ and $\bar{u}_s > \vep$}{
    $s=s+1$\;\label{apl:while_loop1}
    Call the gap reduction $\Gcal(\xbf^{s-1}$, $\bar{l}_{s-1}$, $\eta$, $\theta)$ and obtain $\xbf^s, \bar{l}_{s}$\;
    Set $\bar{u}_s=v(\xbf^s, \eta)$\;\label{apl:while_loop3}
    }
    \KwOut{($\xbf^s, \bar{l}_s$)}
\end{algorithm}
To satisfy \eqref{eq:subprob-relative-error}, 
we develop the APL method (Algorithm~\ref{alg:apl}), which repeatedly calls Algorithm~\ref{alg:gap-reduce} to adjust the bounds $v_k^{\rm L}$ and  $v_k^{\rm U}$.  APL will be used many times in our subsequent development as a subroutine $\Acal(\xbf^0, \bar{l}_0, \eta, \theta, \alpha)$.
It should be noted that our implementation has some differences from the original APL. We terminate the algorithm when either the relative optimality gap~\eqref{eq:subprob-relative-error} or the estimated upper bound of $V(\eta)$ falls below the target threshold. 
Additionally, instead of computing the initial lower bound $\bar{l}_0$  within APL, we provide it as an input parameter to allow a warm start. 
The complexity of Algorithm~\ref{alg:apl} is derived as follows. 
\begin{theorem}\label{thm:complexity-APL}
The total number of iterations of Algorithm~\ref{alg:apl} is bounded by 
    $$S= \max\left\{0, \left\lceil\log_{1/\gamma} \Brbra{\frac{2\alpha(\bar{u}_0-\bar{l}_0)}{(\alpha-1)}\min\bcbra{\frac{1}{V(\eta)},\frac{1}{\vep}}}\right\rceil\right\}$$,
where $\gamma = {(1+\theta)}/{2}$, and $(\blue{\bar{u}_0},\blue{\bar{l}_0})$ are given in Algorithm~\ref{alg:apl}.  Additionally, the number of calls to APL gap reduction is bounded by 
\[S + \frac{1}{1-\gamma^{2/(1+3\rho)}}\Big(\frac{2^{\rho+\blue{2}}3^{(1-\rho)/2}\bar{M}D_\Xcal^{(\rho+1)/2}  \alpha}{(1+\rho)\theta(\alpha-1)}\min\bcbra{\frac{1}{V(\eta)},\frac{1}{\vep}}\Big)^{{2}/{(1+3\rho)}}.\]
\end{theorem}

\subsection{APL-based Root-finding}\label{subsec:practical-root-finding}

In this subsection, we develop practical level-set methods for solving problem~\eqref{pb:func-constraint}, which repeatedly use APL for solving the convex subproblem.

The first step is to identify an appropriate initial value, $\eta_0 < \blue{f^*}$. To achieve this, we introduce an initialization phase~(Algorithm~\ref{alg:phase-I} in Appendix~\ref{subsec:initialization}).
In essence, this routine will either identify the problem as degenerate, in which case a near-optimal solution $\xbf^0$ is found without the need to search for $\blue{f^*}$, or it will return $(\xbf^0, \eta_0, l_0)$  such that  $\eta_0 < \blue{f^*}$ and \eqref{eq:subprob-relative-error} holds with $u_0 = v(\xbf^0, \eta_0)$. Further details are provided in Appendix. 
For the remainder of this discussion, we assume that the problem is non-degenerate and that the required values $\xbf^0$, $\eta_0$, and $l_0 \ge 0$ have been obtained.

\begin{algorithm}[h]
    \small
    \KwIn{$\xbf^0$, $\eta_0<\blue{f^*}$, $l_0$, $\alpha>1$, $\beta\in(0,1)$, $\gamma\in(1/2,1)$,
    $\vep\in(0,\infty)$;}
    Set $u_0=v(\xbf^0, \eta_0), \blue{l_{-1}=l_{0},u_{-1}=2l_{-1}}$ and $t=0$\;
    \While{$u_t>\vep$}{
    Set $t=t+1$,  $\eta_{t}=\eta_{t-1}+\beta l_{t-1}$\label{alg:line:update_etat}\;
    Set {$\tilde{l}_{t}=\max\{1-\beta, \blue{1 + \frac{l_{t-1}-u_{t-2}}{l_{t-2}}}\}\cdot l_{t-1}$},
    $\tilde{u}_t=v(\xbf^{t-1}, \eta_t)$\label{alg:line:update_l_tilde}\;
    Compute $(\xbf^{t}, l_t)=\Acal(\xbf^{t-1},\tilde{l}_{t}, \eta_t, 2\gamma-1, \alpha)$, set $u_t=v(\xbf^t, \eta_t)$\;
    }
    \KwOut{$\xbf^{t}$;}
    \caption{APL-based {\ifix{}}}\label{alg:prac-ifp}
    \end{algorithm}
    \begin{algorithm}[h!]
        \small
        \KwIn{$\xbf^0$, $\eta_0<\blue{f^*}$, $l_0$, \blue{$\beta\in(1/2,1]$, $\alpha\in(1,2\sqrt{\beta})$, $\gamma\in(1/2,1)$}, $\vep\in(0,\infty)$;}
        Set $u_0=v(\xbf^0, \eta_0)$ and $t=0$\;
        \While{$u_t > \vep $}{
        Set $t=t+1$\;
        \lIf{$t=1$}{Set $\eta_t = \eta_{t-1}+\beta l_{t-1}$}
        \lElse{Set $\eta_t=\eta_{t-1} +\beta \cdot \max\bcbra{1, - \frac{\eta_{t-2}-\eta_{t-1}}{u_{t-2}-l_{t-1}}}l_{t-1}$\label{alg:line:update_eta_secant}}
        Set $\tilde{l}_{t}=(1-\beta)l_{t-1}$, $\tilde{u}_t=v(\xbf^{t-1}, \eta_t)$\label{alg:line:update_l_tilde_secant}\;
        Compute $(\xbf^t,l_t)=\mcal A(\xbf^{t-1},\tilde{l}_t,\eta_t,2\gamma - 1, \alpha)$ and set $u_t=v(\xbf^{t}, \eta_t)$\;
        }
        \KwOut{$\xbf^{t}$}
        \caption{APL-based Truncated {\isecant}}\label{alg:secant-2}
        \end{algorithm}

\blue{In each call to the APL method, we need to estimate the lower and upper bounds of $V(\eta_t)$. We use the previous iterate $\xbf^{t-1}$ to compute an upper bound: $\tilde{u}_t=v(\xbf^{t-1}, \eta_t)$ and exploit the convexity of $V(\cdot)$ to obtain the lower bound $\tilde{l}_t$.}
We now present the fixed-point iteration in Algorithm~\ref{alg:prac-ifp} and the APL-based inexact secant method in Algorithm~\ref{alg:secant-2}, both of which solve the subproblem approximately using the APL method described above. 
The convergence properties of \blue{Algorithm~\ref{alg:prac-ifp} and Algorithm~\ref{alg:secant-2}} are summarized in the following theorem.
\begin{theorem}\label{thm:fixed_point_secant_all_complexity}
\blue{Suppose the parameters in Algorithm~\ref{alg:prac-ifp} and Algorithm~\ref{alg:secant-2} satisfy their respective requirements. Then the following results hold:}
\begin{enumerate}
\item \blue{In both algorithms, $\tilde{l}_t$ provides a valid nonnegative lower bound on $V(\eta_t)$ for $t>0$.}
\item \blue{When either algorithm terminates at iteration $T$ with $u_T\leq \vep$, letting $\bar{V}^{\prime}=\min\bcbra{\xi:\xi\in \partial V(f^*)}$, the solution satisfies: $f^*-\eta_T\leq \frac{\vep}{-\bar{V}^{\prime}}$, $f(\xbf^{T})-f^*\le\vep$, and $\norm{[\gbf(\xbf^{T})]_{+}}_\infty\le\vep$.} 
\item Let $\ybf^*\in\Rbb_+^m$ be a vector of Lagrange multipliers at an optimal solution. For Algorithm~\ref{alg:prac-ifp}, the total number of gap reduction calls in the while loop is bounded by $\Ocal\big((\norm{\ybf^{*}}_{1}+1)\log(\norm{\ybf^{*}}_{1}+1)\,(1/\vep)^{2/(1+3\rho)}\big)$. For Algorithm~\ref{alg:secant-2}, the total number of gap reduction calls in the while loop is bounded by $\Ocal\big(\min\{(\norm{\ybf^*}_1+1)\log(\norm{\ybf^*}_1+1), \log(1/\vep)\} \cdot (1/\vep)^{2/(1+3\rho)}\big)$.
\end{enumerate}
\end{theorem}

\begin{remark}
    In view of Theorem~\ref{thm:fixed_point_secant_all_complexity}, the APL-based secant method exhibits two distinct performance behaviors. 
    For well-conditioned problems, where $\norm{\ybf^*}_1\log\norm{\ybf^*}_1=\Ocal(\log(1/\vep))$,  the algorithm performs comparably to the APL-based fixed-point iteration. 
    On the other hand, for ill-conditioned problems with significantly large  $\norm{\ybf^*}$, the robustness of the secant method ensures that it still achieves a near-optimal complexity bound of 
    \blue{${\Ocal}\brbra{(1/\vep)^{2/(1+3\rho)}\log(1/\vep)}$}. 
\end{remark}

\section{Numerical Study}\label{sec:numerical}
In this section, we evaluate the performance of our optimization algorithms across various tasks. \blue{The associated codes are available in~\cite{Qi2026}}. First, we assess the advantage of our acceleration strategy in \apm{} by comparing it to the non-accelerated version~\cite{devanathan2023polyak}. Sections~\ref{sec:socp} and~\ref{sec:LMI} focus on the Karush-Kuhn-Tucker (KKT) systems for Second-Order Cone Programming (SOCP) and Linear Matrix Inequalities (LMI), respectively. Solving the KKT system is equivalent to addressing a penalty problem where the optimal value  $f^*$  is zero.
Second, we explore more general convex optimization problems in the subsequent sections. Section~\ref{sec:convex_qcqp} examines a standard convex Quadratically Constrained Quadratic Programming (QCQP) problem, while Section~\ref{sec:neyman_pearson} investigates the classical Neyman-Pearson classification problem. Extensive experiments demonstrate the superiority of our algorithms. All experiments were conducted on a Mac mini M2 Pro with 32GB of RAM.

\subsection{\label{sec:socp}SOCPs}
We consider the following conic program~\citep{devanathan2023polyak} and its dual problem:
\begin{equation}
    \begin{aligned}
        \min\  & c^{\top}u \\
        \text{s.t.}\  & Au = b,\, u \in \mathcal{K},
    \end{aligned}
    \qquad
    \begin{aligned}
        \max\  & b^{\top}v \\
        \text{s.t.}\  & c - A^{\top}v = s,\, s \in \mathcal{K}^{*},
    \end{aligned}
\end{equation}
\blue{where $\mcal K$ denotes a convex cone and $\Kcal^*$ its dual cone.} This problem can be reformulated as the following penalty problem:
\begin{equation}\label{eq:kkt_socp}
\begin{aligned}
    \min_{\xbf} \  d_{\mcal K}(\ubf) + d_{\mcal K^*}(\sbf)
    \  \st \xbf=\left[\begin{array}{ccc}
\ubf^{\top} & \vbf^{\top} & 1\end{array}\right]^{\top}, \dbf=\left[\begin{array}{ccc}
\sbf^{\top} & 0 & 0\end{array}\right]^{\top},
\ \  E=\left[\begin{array}{ccc}
0 & -A^{\top} & \cbf\\
A & 0 & -\bbf\\
-\cbf^{\top} & \bbf^{\top} & 0
\end{array}\right],\ E\xbf = \dbf,
\end{aligned}
\end{equation}
\blue{where $d_{\mcal K}(\ubf)$ denotes the distance from $u$ to $\mcal K$, and $d_{\mcal K^{*}}(s)$ is the distance from $s$ to $\mcal K^*$. 
In the case of SOCP, $\Kcal$ is the second-order cone, which is self-dual, i.e., $\Kcal^* = \Kcal$.} If the original SOCP admits an optimal solution, then the optimal objective value of~\eqref{eq:kkt_socp} is zero. 

The data $(A, \bbf, \cbf)$ are generated in a manner similar to~\cite{devanathan2023polyak}. Specifically, we first sample a vector $z$ from the standard normal distribution and project it onto the second-order cone $\mcal K$ to obtain $\ubf$. We then set $\sbf = \ubf - \zbf$, which ensures that $\sbf$ lies in the dual cone $\mcal K^*$. Next, both $\vbf$ and $A$ are sampled from the standard normal distribution. Finally, we set $\bbf = A\ubf$ and $\cbf = \sbf + A^\top \vbf$, thus completing the data generation process.

We compare r\apm{} (Algorithm~\ref{alg:restartedAcc}) with {\PMM}~\citep{devanathan2023polyak} by using two problems of different scales. 
We set $\theta=1/2$ in r\apm{}. Both {r\APM} and  {\PMM} require solving a quadratic subproblem formulated as~\eqref{eq:lcqp} at each iteration, which is then solved by the commercial solver Mosek~\citep{aps2019Mosek}. In Figure~\ref{fig:socp_exp}, we plot the convergence of the compared algorithms with respect to the number of gradient evaluations. The red line represents {\PMM}, while the black line represents {r\APM}. We evaluated the performance of the two algorithms under varying bundle sizes ($B$), corresponding to $k_0$ in a limited-memory setup for selecting $X_k$.
For small bundle sizes (i.e., $B\le 5$), we note that {r\apm{}} exhibits much faster convergence than \PMM{}, highlighting the advantages of using momentum acceleration. To further explore the effect of bundle size, we increase it to 100 to examine its impact on convergence rates. It is important to note that using such a large bundle size in practice may be impractical, as the resulting subproblem becomes highly constrained.
As the bundle size increases, both algorithms exhibit faster convergence rates, as expected.  
\blue{We also report the runtime comparison in Table~\ref{tab:time_comparison}. Although r{\APM} introduces additional operations for acceleration, its runtime remains comparable to that of $\PMM$ under the same iteration budget.} 


\begin{figure}
    \begin{centering}
        \begin{minipage}[t]{0.5\columnwidth}
            \includegraphics[width=7.78cm]{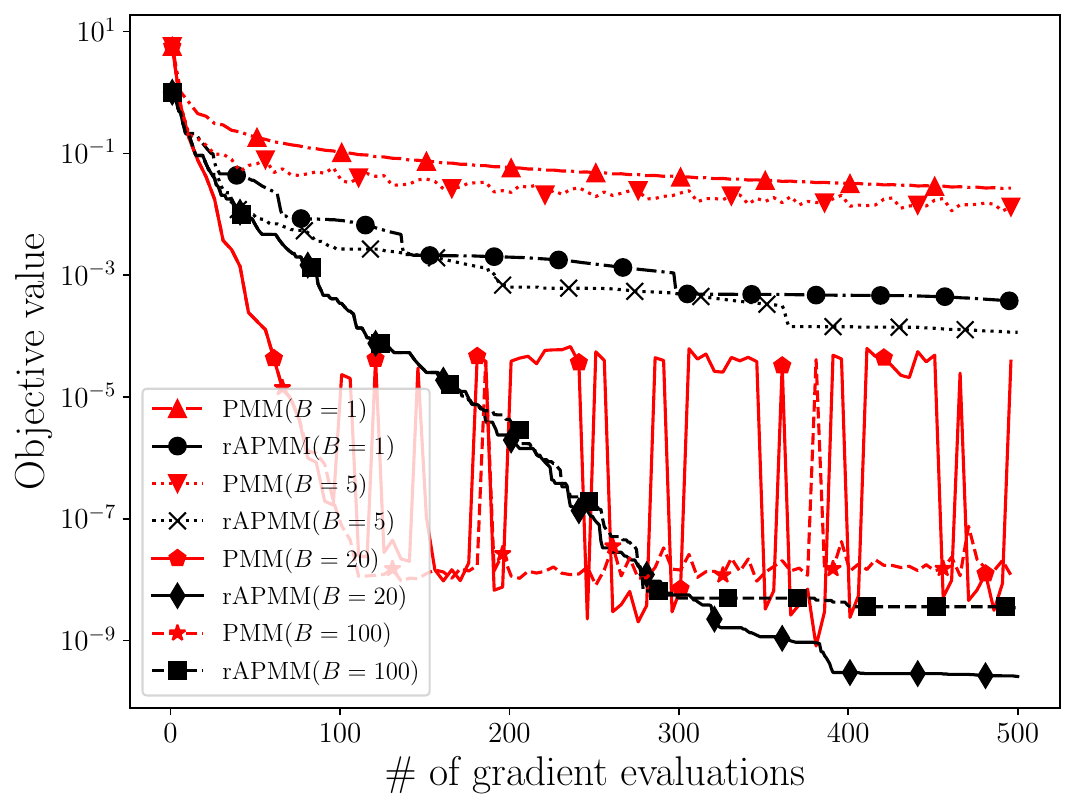}
        \end{minipage}
        \begin{minipage}[t]{0.5\columnwidth}
            \includegraphics[width=8cm]{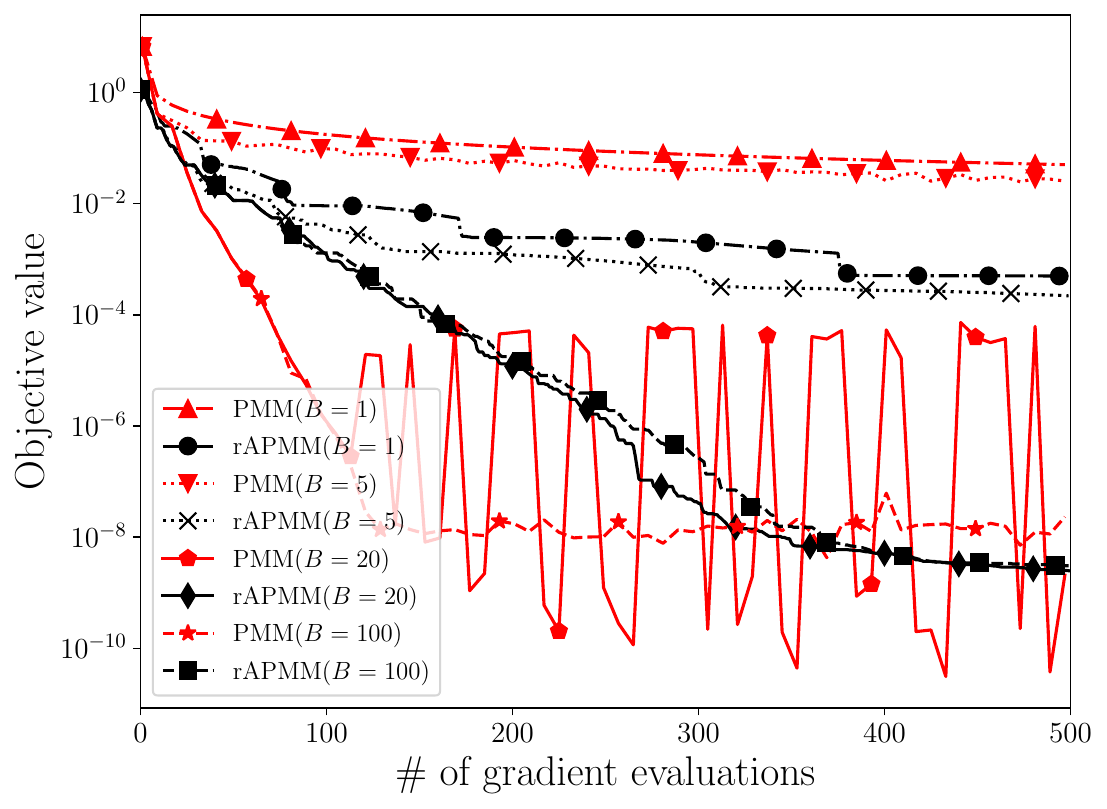}
        \end{minipage}
    \end{centering}
    \caption{\label{fig:socp_exp} SOCP convergence. Left: 500 variables, 200 equality constraints, and 10 cones, each of dimension 50. Right: 1000 variables, 800 equality constraints, and 10 cones, each of dimension 100.}
\end{figure}


\begin{table}[htbp]
  \centering
  \small
  \caption{\blue{\label{tab:time_comparison}Time (seconds) comparison of experiments shown in Figure~\ref{fig:socp_exp} and~\ref{fig:lmi_exp}}}
    \begin{tabular}{lrrrr|rrrr}
    \toprule
    \multicolumn{1}{c}{\multirow{2}[4]{*}{Bundle Size}} & 1     & 5     & 20    & 100   & 1     & 5     & 20    & 100 \\
\cmidrule{2-9}          & \multicolumn{4}{c|}{SOCP Case 1} & \multicolumn{4}{c}{SOCP Case 2} \\
    \midrule
    PMM   & 54.18 & 77.39 & 86.84 & 102.31 & 584.14 & 495.67 & 546.86 & 848.18 \\
    rAPMM & 64.28 & 78.01 & 90.88 & 96.53 & 483.73 & 561.79 & 590.76 & 958.44 \\
\cmidrule{2-9}          & \multicolumn{4}{c|}{LMI Case 1} & \multicolumn{4}{c}{LMI Case 2} \\
\cmidrule{2-9}    PMM   & 3.22  & 13.42 & 52.64 & 273.27 & 24.68 & 131.08 & 567.85 & 2791.29 \\
    rAPMM & 3.17  & 13.44 & 52.94 & 239.59 & 24.25 & 131.58 & 513.76 & 2344.38 \\
    \bottomrule
    \end{tabular}%
  \label{tab:addlabel}%
\end{table}%

\subsection{\label{sec:LMI}LMIs}
The second experiment focuses on a Linear Matrix Inequality (LMI) problem, which seeks a symmetric matrix  $X\in \mbb {S}^{q \times q}$ such that the following inequalities hold:
\begin{equation}\label{eq:feasibility_lmi}
    X\succeq I, \ \ A_i^\top X + XA_i \preceq 0,\ \  i = 1,\cdots, k,
\end{equation}
where the inequalities are with respect to the positive semidefinite cone, and $A_i\in \mbb R^{q\times q}$. This can also be expressed as the following optimization problem:
\begin{equation}\label{eq:lmi_obj}
    \min_{X\in \mbb S^{q\times q}} 0,\ \  \st \sigma_{\max}(I-X)\le 0, \ \  \sigma_{\max}(A_i^\top X + XA_i)\le 0,\ \ \forall i=1,\ldots,k.
\end{equation}
where  $\sigma_{\max}(\cdot)$  denotes the largest singular value.
If problem~\eqref{eq:feasibility_lmi} has a feasible solution, then the optimal objective value of the above problem~\eqref{eq:lmi_obj} is 0. The data generation method for $(A_i,i=1,\cdots, k)$ follows the approach proposed by~\cite{devanathan2023polyak}. Specifically, we first sample the $q\times q$ matrices $B_i$ and $C_i$ from standard normal distributions. Then each $A_i = F^{-1}\brbra{-B_i B_i^\top +C_i-C_i^\top}F$, which ensures that $X=F^\top F$ satisfies $A_i^\top X+XA_i \preceq 0,i=1,\cdots, k$. Since $X\succeq 0$, we scale it to satisfy $X\succeq I$, thereby forming a feasible instance of problem~\eqref{eq:feasibility_lmi}. We compare the convergence performance of {r\APM} and {\PMM} under two different data sizes, as shown in Figure~\ref{fig:lmi_exp}. Consistent with the finding in Section~\ref{sec:socp}, both algorithms exhibit faster convergence as the bundle size increases. As shown in the right part of Figure~\ref{fig:lmi_exp}, {r\APM} not only converges faster than {\PMM} but also exhibits a more stable convergence trajectory.
\begin{figure}
    \begin{centering}
        \begin{minipage}[t]{0.5\columnwidth}
            \includegraphics[width=8cm]{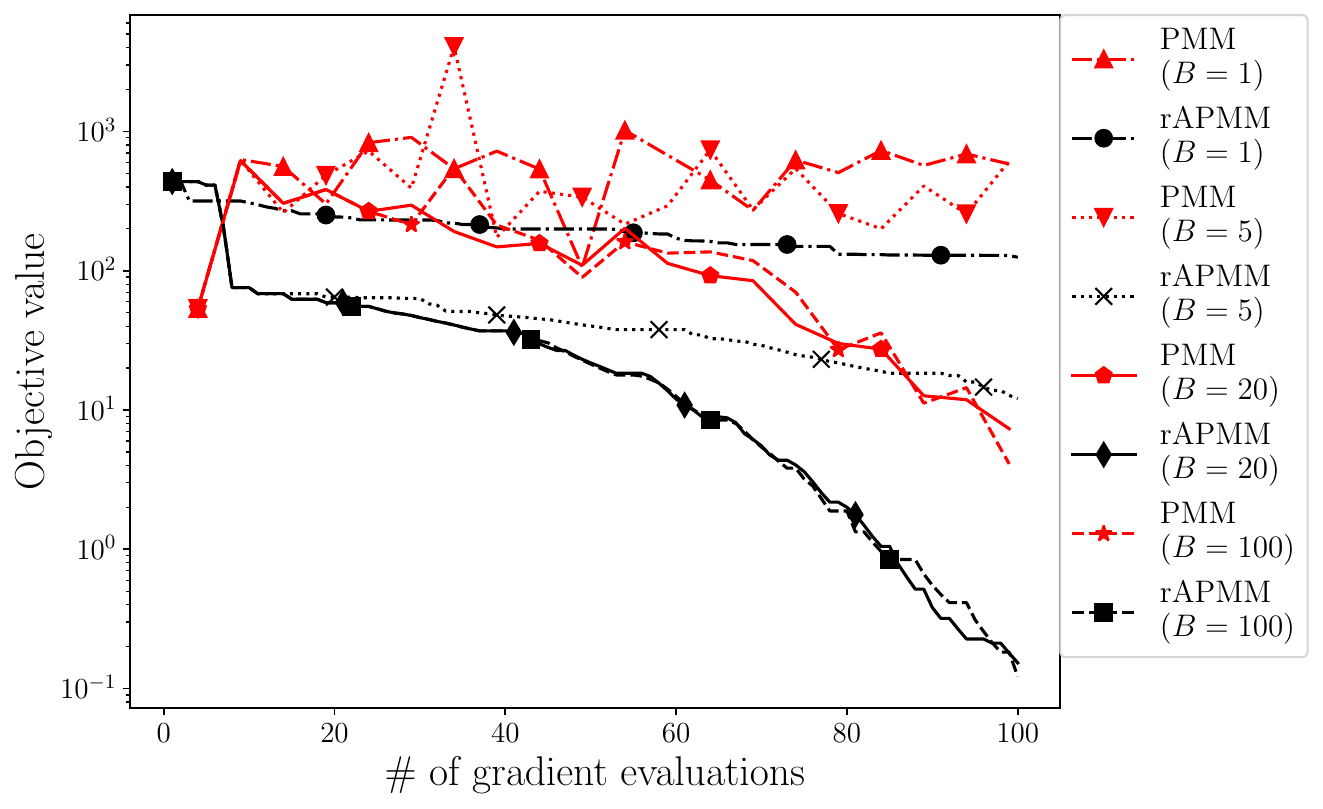}
        \end{minipage}
        \begin{minipage}[t]{0.5\columnwidth}
            \includegraphics[width=8cm]{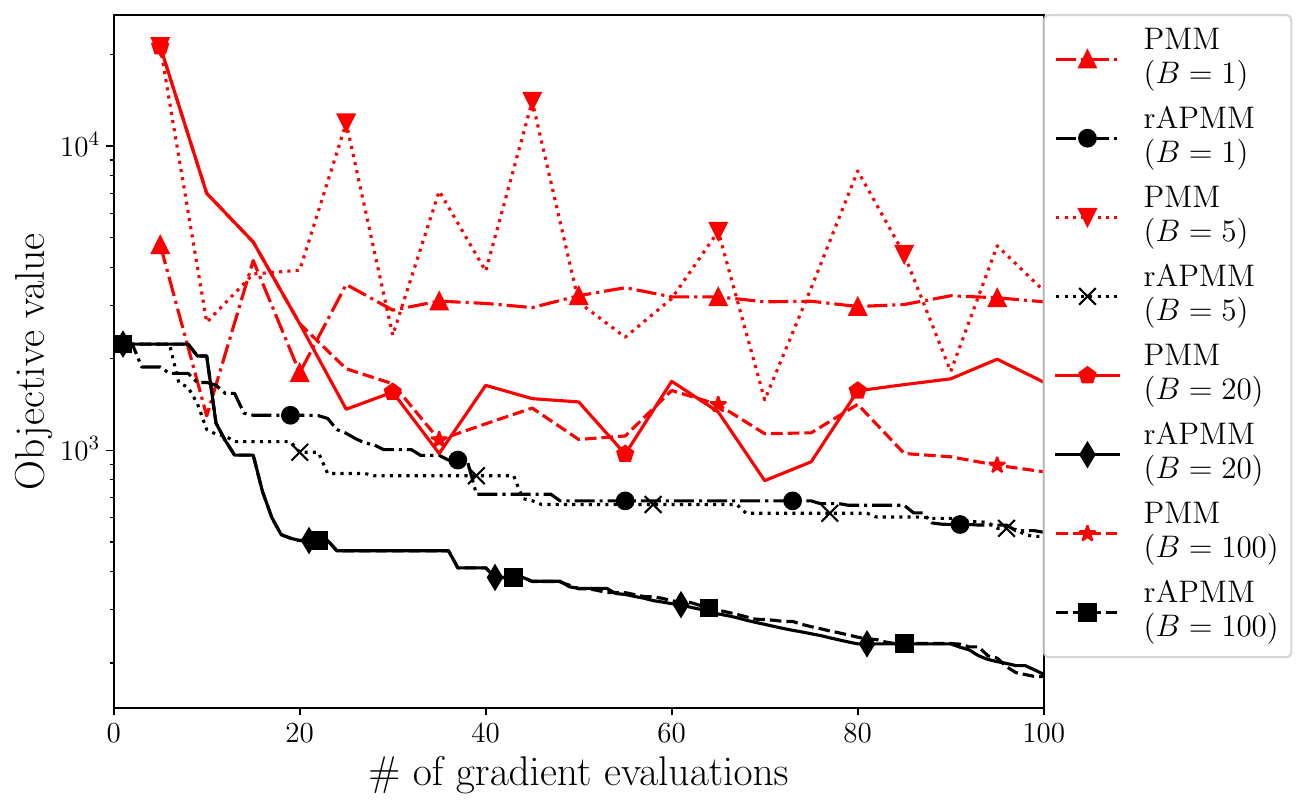}
        \end{minipage}
    \end{centering}
    \caption{\label{fig:lmi_exp}
    Convergence performance on LMI.
    Left:  $(q,k) = (20,10)$; Right:  $(q,k) = (40,20)$. }
\end{figure}

\subsection{Convex QCQPs\label{sec:convex_qcqp}}
\begin{figure}
    \begin{centering}
        \begin{minipage}[t]{0.5\columnwidth}
            \includegraphics[width=8cm]{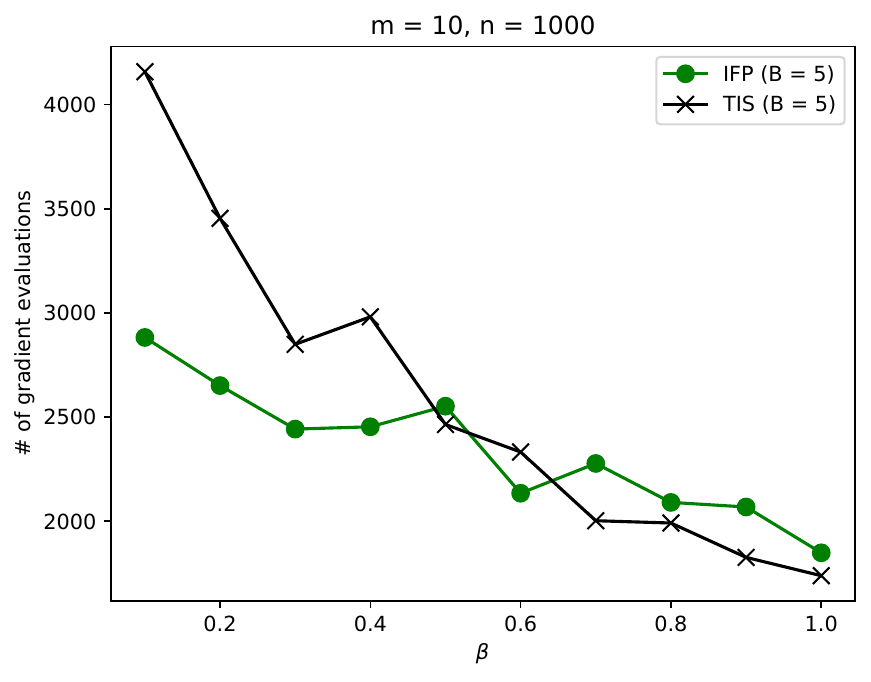}
        \end{minipage}
        \begin{minipage}[t]{0.5\columnwidth}
            \includegraphics[width=8cm]{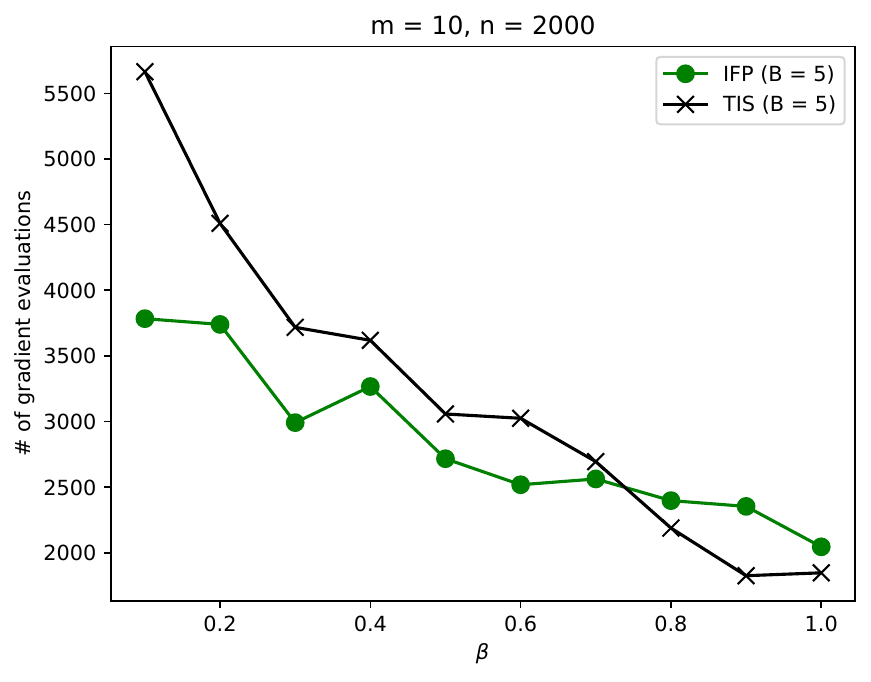}
        \end{minipage}
    \end{centering}
    \caption{\label{fig:iteration_beta}Gradient evaluations vs.\ $\beta$ for convex QCQP ($\vep = 10^{-3}$). Left: $(m,n)=(10,1000)$; right: $(10,2000)$. IFP: Inexact Fixed Point; TIS: Truncated Inexact Secant.}
\end{figure}
We consider the following convex QCQP problem:
\begin{equation}
\min_{\xbf \in \mbb R^n, \rm{ub}\succeq \xbf \succeq \rm{lb}}  \   \frac{1}{2} \xbf^\top Q_0 \xbf + \cbf^\top_0 \xbf \ \ \st\ \ \ \frac{1}{2} \xbf^\top Q_i \xbf + \cbf^\top_i\xbf + d_i \leq 0,\ i = 1, \ldots, m,
\end{equation}
where every $Q_i\succeq 0,i=0,\ldots,m$, each $\cbf_i, i=0,\ldots,m$ is sampled from Gaussian distribution, $d_i$ is set as constant $10$ and set $\rm{lb}=-10,\rm{ub}=10$.
We assume the optimal value to this problem is unknown in advance. Our goal is to evaluate the performance of the APL-based {\ifix{}{}} (Algorithm~\ref{alg:prac-ifp}) and the APL-based Truncated {\isecant} (Algorithm~\ref{alg:secant-2}).

We first investigate how the stepsize $\beta$ influences the performance of the algorithms. 
There is a potential trade-off between optimizing the efficiency of the outer loop and the inner loop. A larger $\beta$  leads to a more aggressive update of $ \eta_t $, but it may result in a less accurate initial lower bound  $\tilde{l}_t$  for APL. On the other hand, a smaller  $\beta$  tends to provide a better initial estimate of the lower bound, potentially improving the accuracy of the initial gap estimation for APL.
We use Figure~\ref{fig:iteration_beta} to illustrate the impact of different values of $\beta$ on the number of iterations. In all experiments shown in Figure~\ref{fig:iteration_beta}, we used a limited-memory set with a bundle size of 5 for $X_k$, along with parameters $\gamma = 0.9$ and $\alpha = 1.36$.  Overall, we observe that as $\beta$ increases, the number of iterations tends to decrease.  This suggests that the trade-off favors using a larger stepsize. APL is robust to the initial solution, and setting  $\beta = 1$ is empirically a satisfactory choice.

Next, we compared our two methods against Mosek on large-scale convex QCQPs. We fixed  $m = 10$  and varied the values of $ n$. For each random instance, we ran the algorithms on five different random initial solutions and summarized the results in Table~\ref{tab:convex_qcqp_res}. While our methods involve repeatedly solving subproblems, we observed that these subproblems are solved significantly faster than the original QCQP, providing a substantial advantage in large-scale settings. As shown in Table~\ref{tab:convex_qcqp_res}, when the problem size exceeds 6,000, Mosek fails to find a solution, whereas our methods maintain stable performance with no significant increase in overall computation time. Among the three methods, the inexact secant method demonstrates the best performance. To further evaluate scalability with respect to the number of constraints, we additionally fixed the variable dimension $n=4000$ and varied $m$ from $20$ to $100$, with results reported in Table~\ref{tab:qcqp_exp_diff_m}. The results show that Mosek becomes unavailable once $m\ge 30$, whereas both IFP and TIS remain robust and continue to deliver high-quality solutions; moreover, TIS is consistently faster than IFP, further demonstrating the practical advantage of our approach in higher-constraint regimes.

\begin{table}[!t]
  \centering
  \small
  \caption{\label{tab:convex_qcqp_res}Convex QCQP results for varying $n$ with fixing constraint number $m=10$ (first row: mean; second row: standard deviation).}
    \begin{tabular}{ccccccc}
    \toprule
\multirow{2}{*}{$n$} & \multicolumn{2}{c}{Mosek} & \multicolumn{2}{c}{IFP} & \multicolumn{2}{c}{TIS}\tabularnewline
\cmidrule(lr{8pt}){2-7}
 & optVal & time & optVal & time & optVal & time \tabularnewline
    \midrule
    \multirow{2}[2]{*}{4000} & -3.64E+01 & 2.40E+02 & -3.64E+01 & 1.13E+02 & -3.64E+01 & 7.67E+01 \\
          & $\pm$0.00E+00 & $\pm$1.21E+01 & $\pm$1.46E-04 & $\pm$1.85E+00 & $\pm$3.70E-04 & $\pm$3.25E+00 \\
    \midrule
    \multirow{2}[2]{*}{5000} & -4.25E+01 & 4.28E+02 & -4.25E+01 & 1.43E+02 & -4.25E+01 & 9.55E+01 \\
          & $\pm$0.00E+00 & $\pm$1.59E+01 & $\pm$1.62E-04 & $\pm$6.05E+00 & $\pm$1.89E-04 & $\pm$4.24E+00 \\
    \midrule
    \multirow{2}[2]{*}{6000} & -4.67E+01 & 7.68E+02 & -4.67E+01 & 1.68E+02 & -4.67E+01 & 1.10E+02 \\
          & $\pm$0.00E+00 & $\pm$2.19E+01 & $\pm$1.71E-04 & $\pm$4.75E+00 & $\pm$4.23E-04 & $\pm$3.41E+00 \\
    \midrule
    \multirow{2}[2]{*}{7000} &   -    &   -    & -5.23E+01 & 1.86E+02 & -5.23E+01 & 1.20E+02 \\
          &   -    &   -    & $\pm$1.03E-04 & $\pm$3.37E+00 & $\pm$2.93E-04 & $\pm$4.24E+00 \\
    \bottomrule
    \end{tabular}%
\end{table}

\begin{table}[!t]
    \centering
    \small
    \caption{Convex QCQP results for varying $m$ with fixing variable dimension $n=4,000$ (first row: mean; second row: standard deviation).}
      \begin{tabular}{cllrrrr}
      \toprule
      \multirow{2}[4]{*}{$m$} & \multicolumn{2}{c}{Mosek} & \multicolumn{2}{c}{IFP} & \multicolumn{2}{c}{TIS} \\
  \cmidrule{2-7}          & optVal & time  & \multicolumn{1}{l}{optVal} & \multicolumn{1}{l}{time} & \multicolumn{1}{l}{optVal} & \multicolumn{1}{l}{time} \\
      \midrule
      \multirow{2}[2]{*}{20} & \multicolumn{1}{r}{-3.31E+01} & \multicolumn{1}{r}{5.50E+02} & -3.31E+01 & 1.95E+02 & -3.31E+01 & 1.84E+02 \\
            & \multicolumn{1}{r}{$\pm$3.82E-01} & \multicolumn{1}{r}{$\pm$1.17E+01} & $\pm$3.82E-01 & $\pm$8.18E+00 & $\pm$3.82E-01 & $\pm$8.42E+00 \\
      \midrule
      \multirow{2}[2]{*}{30} & - & - & -3.26E+01 & 2.43E+02 & -3.26E+01 & 2.10E+02 \\
            & - & - & $\pm$4.64E-01 & $\pm$2.29E+01 & $\pm$4.64E-01 & $\pm$1.22E+01 \\
      \midrule
      \multirow{2}[2]{*}{60} & - & - & -3.19E+01 & 4.79E+02 & -3.19E+01 & 4.25E+02 \\
            & - & - & $\pm$3.75E-01 & $\pm$4.11E+01 & $\pm$3.75E-01 & $\pm$2.85E+01 \\
      \midrule
      \multirow{2}[2]{*}{80} & - & - & -3.17E+01 & 6.40E+02 & -3.17E+01 & 5.53E+02 \\
            & - & - & $\pm$3.76E-01 & $\pm$7.54E+00 & $\pm$3.77E-01 & $\pm$2.05E+01 \\
      \midrule
      \multirow{2}[2]{*}{100} & - & - & -3.16E+01 & 8.75E+02 & -3.16E+01 & 7.79E+02 \\
            & - & - & $\pm$3.89E-01 & $\pm$6.75E+01 & $\pm$3.89E-01 & $\pm$5.67E+01 \\
      \bottomrule
      \end{tabular}%
    \label{tab:qcqp_exp_diff_m}%
  \end{table}%


\subsection{Neyman-Pearson classification\label{sec:neyman_pearson}}
The final experiment involves the Neyman-Pearson classification (NPC) problem~\citep{tong2016survey}. In binary classification, this problem focuses on minimizing the classification loss of one class while controlling the classification loss of the other. Using the logistic loss as a surrogate for classification loss, the binary NPC problem can be formulated as follows:
\begin{equation}\label{eq:binary_np}
    \begin{aligned}
        \min_{\wbf \in \mcal W(r)}&\ \  \frac{1}{n_1}\sum_{i=1}^{n}\onebf(y_i=1)\log\brbra{1+\exp(-y_i \cdot \xbf_i^{\top} \wbf)} + \frac{\rho}{2}\norm{\wbf}^2\\
        \st &\ \  \frac{1}{n_{-1}}\sum_{i=1}^{n}\onebf(y_i=-1)\log\brbra{1+\exp(-y_i \cdot \xbf_i^{\top} \wbf)} \leq \kappa,
    \end{aligned}
\end{equation}
where $\mcal W(r):=\{\wbf \mid \norm{\wbf}^2 \leq r^2\}$, $n_{1}$ and $n_{-1}$ are the number of samples in the categories $1$ and $-1$, respectively.

In multi-class classification, we extend this formulation by controlling the classification loss for each class while minimizing the overall classification loss across all classes. Using the cross-entropy loss, the multi-class NPC problem can be formulated as follows~\eqref{eq:multi_np}:
\begin{equation}\label{eq:multi_np}
\begin{aligned}
    \min_{\wbf_j, j=1,\ldots,J}&\ \  -\frac{1}{n}\sum_{i=1}^{n}\sum_{j=1}^{J}\onebf(y_i = j)\log(p_{ij})\ \ \st -\frac{1}{n_j}\sum_{i=1}^{n}\onebf(y_i = j)\log(p_{ij})\leq \kappa_j, \Wbf \in \mcal W(r)
\end{aligned}
\end{equation}
where $n_j = \sum_{i=1}^{n}\onebf(y_i = j)$, $p_{ij}=\frac{\exp(\xbf_i^{\top} \wbf_j)}{\sum_{j=1}^{J}\exp(\xbf_i^{\top} \wbf_j)}$, $\Wbf=[\wbf_1^\top,\ldots,\wbf_J^\top]^{\top}$ and $\mcal W(r):=\{\Wbf \mid \norm{\Wbf}^2\leq r^2\}$. 
\begin{figure}
    \begin{centering}
        \begin{minipage}[t]{0.48\columnwidth}
            \includegraphics[width=7.5cm]{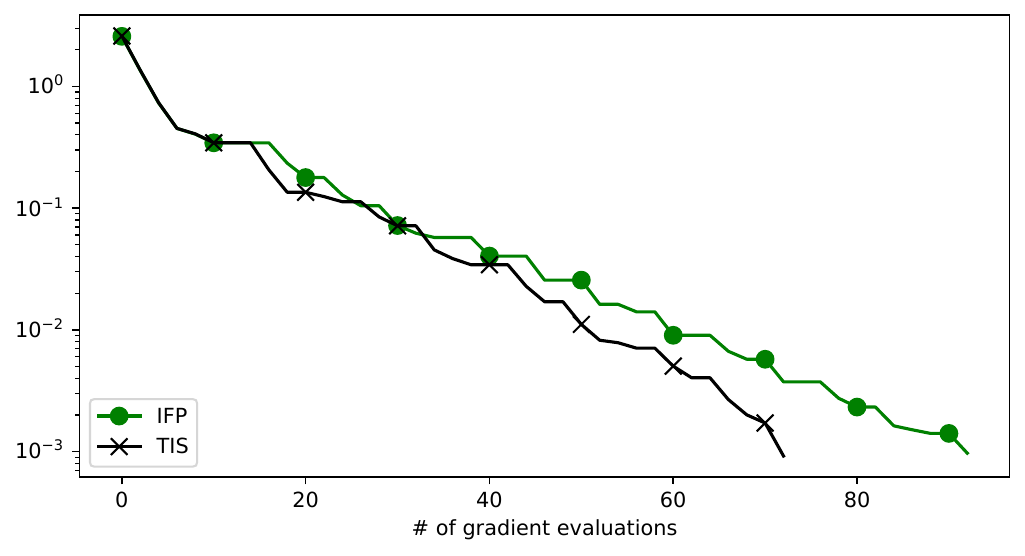}
        \end{minipage}
        \begin{minipage}[t]{0.25\columnwidth}
            \includegraphics[width=4.2cm]{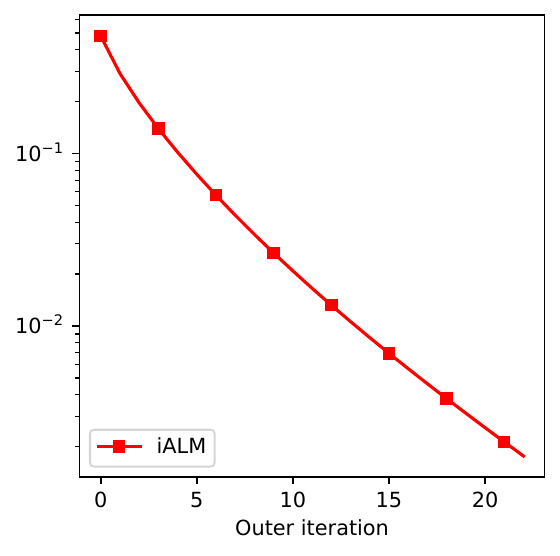}
        \end{minipage}
        \begin{minipage}[t]{0.25\columnwidth}
            \includegraphics[width=4.2cm]{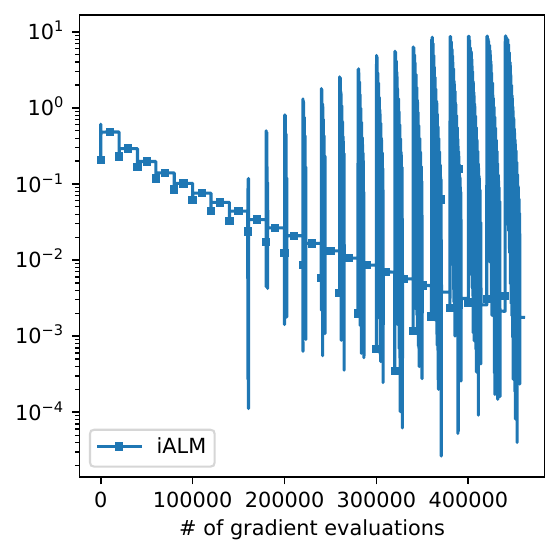}
        \end{minipage}
    \end{centering}
    \caption{Results on NPC. \label{fig:binary_npc}$y$-axis: $\max\{|f(\xbf_k)-f^*|,\{[g_i(\xbf_k)]_{+}\}_i\}$. $x$-axis: first panel uses gradient evaluations to meet the accuracy of $\leq 10^{-3}$ (our methods); the second and third panels report iALM outer iterations and total gradient evaluations, respectively, under a complementarity tolerance of $\leq 10^{-3}$.
    }
\end{figure}
\begin{figure}[h]
    \centering
    \includegraphics[width=0.5\linewidth]{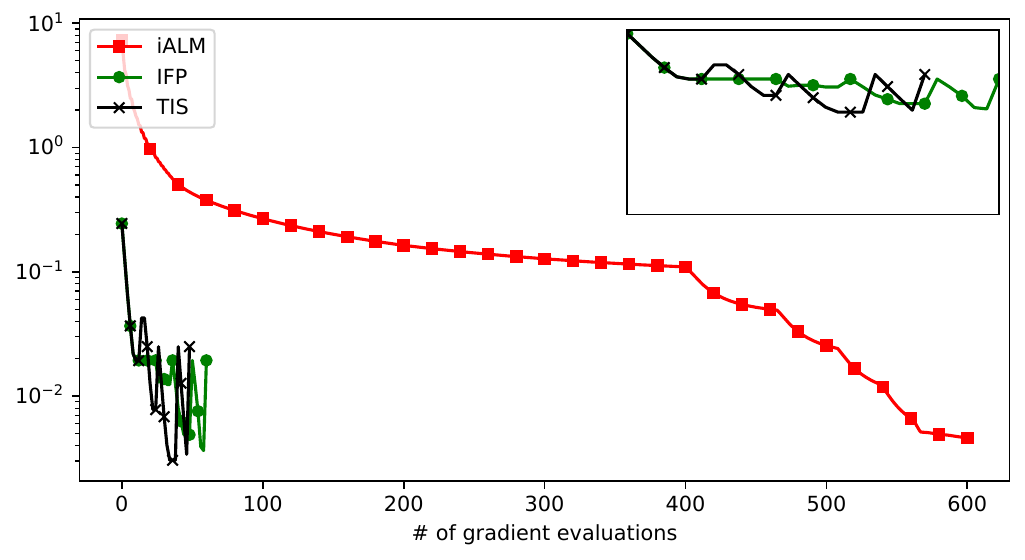}
    \caption{Results on multi-class NPC. $y$-axis: $\max\{|f(\xbf_k)-f^*|,\{[g_i(\xbf_k)]_{+}\}_i\}$. $x$-axis: gradient evaluations for our methods to reach a tolerance of $\leq 10^{-3}$ and iALM to reach a complementarity tolerance of $\leq 10^{-3}$.}
    \label{fig:multi-npc}
\end{figure}

We compare Algorithm~\ref{alg:prac-ifp}, Algorithm~\ref{alg:secant-2} and the inexact Augmented Lagrangian Method (iALM, \cite{xu2021iteration}). Note that iALM is a double-loop method in which the penalty subproblem is solved inexactly using the accelerated gradient method. For binary NPC, we use the Arcene  dataset~\citep{misc_arcene_167}, setting $\rho = 0.01,\kappa=0.5$ and $r=7$. For multi-class NPC, we use the Dry Bean dataset~\citep{Koklu2020MulticlassCO}, setting  $r=7$  and $\kappa_j = 0.8$, $j=1,2,\ldots, J$. For the binary NPC problem, we observe that iALM requires significantly more calls to the oracles than our methods.  To illustrate this, Figure~\ref{fig:binary_npc} presents three graphs comparing the complexities. While the number of outer-loop iterations for iALM remains below 30, the total iteration number, including the inner loop, is significantly higher.  From the last plot in Figure~\ref{fig:binary_npc}, we observe that iALM exhibits significant oscillations in solving the subproblems. A possible explanation is that the performance of the accelerated gradient method is sensitive to the conditioning of the iALM subproblem, which dynamically changes after each update of the dual variables. This makes tuning the parameters of iALM substantially more challenging. 
In contrast, our algorithms exhibit more stable descent patterns,  and we observe that the secant algorithm demonstrates superior convergence compared to the fixed point iteration. In the experiment of multi-class NPC (shown in Figure~\ref{fig:multi-npc}), we also observe that our APL-based level set methods perform better than iALM, while the performance of the secant method and fixed point iteration are close.

\section{Conclusions}
This paper presents uniformly optimal and parameter-free first-order methods for convex and function-constrained optimization. For problems with a known optimal value, 
we develop an accelerated algorithm that extends the classic Polyak step-size method and attains optimal complexity. When $f^*$ is unknown, we provide, to our knowledge, the first methods that are both parameter-free and uniformly optimal in the H\"{o}lder smooth setting.
Several directions for future research emerge from this work. Extensions to stochastic optimization and the integration of our level-set algorithms with other parameter-free first-order methods~\cite{nesterov2015universal} constitute promising avenues for further investigation.

\section*{Acknowledgments}
Q. Deng was supported in part by the National Natural Science Foundation of China (12571325, 72394364/72394360) and the Natural Science Foundation of Shanghai (24ZR1421300).

\bibliographystyle{abbrvnat}
\bibliography{ref}

\appendix

\section{Missing Proofs in Section~\ref{sec:apmm}}

\subsection{Proof of Theorem~\ref{thm:rate-APM}}
Before proving Theorem~\ref{thm:rate-APM}, we establish some important properties of the generated sequence $\{\xbf^k\}$.
\begin{proposition}\label{prop:stable}
Let $\xbf^{*}$ be an optimal solution of problem~\eqref{pb:func-constraint}, and assume $\xbf^*\in X_k$, for $k=1,2,\ldots$.
Then, $\xbf^{*}$ is a feasible \blue{point} of the subproblem~\eqref{eq:update-xk}.
Moreover, the sequence generated by Algorithm~\ref{alg:APMM} has the following square summable property:
\begin{equation}
\|\xbf^{*}-\xbf^{K}\|^2+\sum_{k=1}^{K}\|\xbf^{k}-\xbf^{k-1}\|^2\le \|\xbf^{*}-\xbf^{0}\|^2.\label{eq:sum-square-xk}
\end{equation}
\end{proposition}
\begin{proof}
The feasibility of $\xbf^{*}$ immediately follows from the convexity:
$\ell_{f}(\xbf^{*},\zbf^{k})\le f(\xbf^{*})=f^{*}, \ \text{and } \ell_{g_i}(\xbf^*, \zbf^k) \le g_i(\xbf^*)\le 0, i\in [m].$
Applying the optimality condition~\citep[Proposition 2.7.]{lan2020first} to~\eqref{eq:update-xk}, we have $\inner{\xbf^k-\xbf^{k-1}}{\xbf-\xbf^k}\ge 0$ for any $\xbf\in X_{k-1}$ such that $\ell_{f}(\xbf,\zbf^{k})\le f^{*}$.
By simple algebraic manipulation, this can be rewritten as
$\|\xbf^{k}-\xbf^{k-1}\|^2 \le \|\xbf-\xbf^{k-1}\|^2-\|\xbf-\xbf^{k}\|^2,\ \forall\xbf\in X_{k-1},\,\ell_{f}(\xbf,\zbf^{k})\le f^{*}.$
Since we have shown the feasibility of $\xbf^{*}$, placing $\xbf=\xbf^{*}$
in the above relation gives 
\[\|\xbf^{k}-\xbf^{k-1}\|^2\le \|\xbf^{*}-\xbf^{k-1}\|^2-\|\xbf^{*}-\xbf^{k}\|^2,\ k=1,2,\ldots,K.\]
Summing up the above relation over $k=1,2,\ldots,K$, we have the desired result.

\end{proof}

\begin{proof}[Proof of Theorem~\ref{thm:rate-APM}]
First, by telescoping and using the definition of $\ybf^{k}$, we
observe that $\ybf^{k}$ can be expressed as the convex combination
of $\xbf^{0},\xbf^{1},\ldots,\xbf^{k}$. In particular, we can write
$\ybf^{k}=\sum_{s=0}^{k}\beta_{s}\xbf^{s}$ where $\beta_{s}\ge0$
($0\le s\le k$) and $\sum_{s=0}^{k}\beta_{s}=1$. Since $\|\xbf^{*}-\cdot\|^2$ is convex,
applying Jensen's inequality and Proposition~\ref{prop:stable}, we have 
\[\|\xbf^{*}-\ybf^{k}\|^2\le\sum_{s=0}^{k}\beta_{s}\|\xbf^{*}-\xbf^{s}\|^2\le \|\xbf^{*}-\xbf^{0}\|^2.
\]
Similarly, we can show $\|\xbf^{*}-\tilde{\ybf}^{k}\|\le \|\xbf^{*}-\xbf^{0}\|$,
where $\tilde{\ybf}^{k}=\alpha_{k}\xbf^{k}+(1-\alpha_{k})\ybf^{k-1}$.
Using convexity again with the update $\zbf^k$ we have $\|\xbf^{*}-\zbf^{k}\|^2\le(1-\alpha_{k}) \|\xbf^{*}-\ybf^{k-1}\|^2+\alpha_{k} \|\xbf^{*}-\xbf^{k-1}\|^2\le \|\xbf^{*}-\xbf^{0}\|^2.$
Thus, we establish a uniform bound on the curvature of the sequence:
\begin{equation}\label{eq:Mk-bound}
M_{i}(\tilde{\ybf}^{k},\zbf^{k})\le\hat{M},\ k=1,2,\ldots,K.
\end{equation}
Using the smoothness condition, we have
\begin{equation}\label{eq:opt-01}
    \begin{aligned}
f(\tilde{\ybf}^{k}) & \le\ell_{f}(\tilde{\ybf}^{k},\zbf^{k})+\frac{M_{0}(\tilde{\ybf}^{k},\zbf^{k})}{1+\rho}\norm{\tilde{\ybf}^{k}-\zbf^{k}}^{1+\rho} \\
 & =\alpha_{k}\ell_{f}(\xbf^{k},\zbf^{k})+(1-\alpha_{k})\ell_{f}(\ybf^{k-1},\zbf^{k})+\frac{M_{0}(\tilde{\ybf}^{k},\zbf^{k})\alpha_{k}^{1+\rho}}{1+\rho}\norm{\xbf^{k}-\xbf^{k-1}}^{1+\rho}\\
 & \le \alpha_{k}\ell_{f}(\xbf^{k},\zbf^{k})+(1-\alpha_{k})f(\ybf^{k-1})+\frac{M_{0}(\tilde{\ybf}^{k},\zbf^{k})\alpha_{k}^{1+\rho}}{1+\rho}\norm{\xbf^{k}-\xbf^{k-1}}^{1+\rho} \\
 & \le\alpha_{k}f^{*}+(1-\alpha_{k})f(\ybf^{k-1})+\frac{M_{0}(\tilde{\ybf}^{k},\zbf^{k})\alpha_{k}^{1+\rho}}{1+\rho}\norm{\xbf^{k}-\xbf^{k-1}}^{1+\rho}.
 \end{aligned}
\end{equation}
Here, the last inequality uses the fact that $\ell_{f}(\xbf^{k},\zbf^{k})-f^*\le v_\ell(\xbf^k,\zbf^k,f^*)\le 0$.
After rearranging terms, we have $
f(\tilde{\ybf}^{k})-f^{*}\le(1-\alpha_{k})\bsbra{f(\ybf^{k-1})-f^{*}}+\frac{M_{0}(\tilde{\ybf}^{k},\zbf^{k})\alpha_{k}^{1+\rho}}{1+\rho}\norm{\xbf^{k}-\xbf^{k-1}}^{1+\rho}.$
Using a similar argument, we have $
g_{i}(\tilde{\ybf}^{k})\le(1-\alpha_{k})g_i(\ybf^{k-1})+\frac{M_{i}(\tilde{\ybf}^{k},\zbf^{k})\alpha_{k}^{1+\rho}}{1+\rho}\norm{\xbf^{k}-\xbf^{k-1}}^{1+\rho},\ i=1,2,\ldots,m.$
Combining these two inequalities, we have $
v(\ybf^{k},f^{*})\le v(\tilde{\ybf}^{k},f^{*})\le(1-\alpha_{k})v(\ybf^{k-1},f^{*})+\frac{M(\tilde{\ybf}^{k},\zbf^{k})\alpha_{k}^{1+\rho}}{1+\rho}\norm{\xbf^{k}-\xbf^{k-1}}^{1+\rho},$
where $M(\tilde{\ybf}^{k},\zbf^{k})=\max_{0\le i\le m}M_{i}(\tilde{\ybf}^{k},\zbf^{k})$.
Multiplying both sides by $\Gamma_{k}$ \blue{and using the relation $\Gamma_{k-1}=\Gamma_k (1-\alpha_k)$}, we obtain 
\begin{equation*}
    \begin{aligned}
        \Gamma_{k}v(\ybf^{k},f^{*}) & \le\Gamma_{k-1}v(\ybf^{k-1},f^{*})+\frac{M(\tilde{\ybf}^{k},\zbf^{k})\alpha_{k}^{1+\rho}\Gamma_{k}}{1+\rho}\norm{\xbf^{k}-\xbf^{k-1}}^{1+\rho},\ k>1.\\
\Gamma_{1}v(\ybf^{1},f^{*}) & \le\Gamma_{1}(1-\alpha_{1})v(\ybf^{0},f^{*})+\frac{M(\tilde{\ybf}^{1},\zbf^{1})\alpha_{1}^{1+\rho}\Gamma_{1}}{1+\rho}\norm{\xbf^{1}-\xbf^{0}}^{1+\rho},\ k=1.
    \end{aligned}
\end{equation*}
Summing up the above result for $k=1,2,\ldots,K$ and using the bound~(\ref{eq:Mk-bound}),
we obtain $\Gamma_{K}v(\ybf^{K},f^{*})\le\Gamma_{1}(1-\alpha_{1})v(\ybf^{0},f^{*})+\frac{\hat{M}_{}}{1+\rho}\sum_{k=1}^{K}\alpha_{k}^{1+\rho}\Gamma_{k}\norm{\xbf^{k}-\xbf^{k-1}}^{1+\rho}.$
Using H\"{o}lder's inequality (i.e., $\inner{\xbf}{\ybf}\le\norm{\xbf}_{2/(1+\rho)}\cdot\norm{\ybf}_{2/(1-\rho)}$),
it follows that $\Gamma_{K}v(\ybf^{K},f^{*}) \le\Gamma_{1}(1-\alpha_{1})v(\ybf^{0},f^{*})+\frac{\hat{M}}{1+\rho}\norm{\cbf_{K}}_{\frac{2}{1-\rho}}\left(\sum_{k=1}^{K}\norm{\xbf^{k}-\xbf^{k-1}}^{2}\right)^{\frac{\rho+1}{2}}.$
On the other hand, using the relation~(\ref{eq:sum-square-xk}), we have $\sum_{k=1}^{K}\norm{\xbf^{k}-\xbf^{k-1}}^{2}\le \|\xbf^{*}-\xbf^{0}\|^2.$
By putting these two results together, we obtain the desired result.

Next, we derive a more specific convergence rate. Setting $\alpha_{k}=\frac{2}{k+1}$, we have $\Gamma_{k}=\frac{k(k+1)}{2}$.
For $\rho < 1$, we have 
$\norm{\cbf_{K}}_{\frac{2}{1-\rho}} =\Big[\tsum_{k=1}^{K}\left(\frac{2^{\rho}k}{(k+1)^{\rho}}\right)^{\frac{2}{1-\rho}}\Big]^{\frac{1-\rho}{2}}\le2^{\rho}\Big[\tsum_{k=1}^{K}\left((k+1)^{1-\rho}\right)^{\frac{2}{1-\rho}}\Big]^{\frac{1-\rho}{2}} \blue{=}2^{\rho}\Big[\tsum_{k=1}^{K}(k+1)^{2}\Big]^{\frac{1-\rho}{2}}$.

When $\rho=1$, it is easy to verify $\norm{\cbf_K}_\infty=\frac{2K}{K+1}\le 2^\rho$, hence the above inequality holds.
Using the basic fact that $\sum_{s=1}^{k}s^{2}=\frac{k(k+1)(2k+1)}{6}\le\frac{k(k+1)^{2}}{3}$,
it follows that $\Gamma_{K}^{-1}\norm{\cbf_{K}}_{\frac{2}{1-\rho}} \le\Gamma_{K}^{-1}2^{\rho}\left[\frac{(K+1)(K+2)^{2}}{3}\right]^{(1-\rho)/2}
 =\frac{2^{\rho+1}}{3^{(1-\rho)/2}}\frac{(K+2)^{(1-\rho)}}{K(K+1)^{(1+\rho)/2}}
 \blue{\leq\frac{2^{\rho+1}}{3^{(1-\rho)/2}}\frac{\brbra{3K}^{1-\rho}}{K(K+1)^{(1+\rho)/2}}} =\frac{2^{\rho+1}3^{(1-\rho)/2}}{K^{\rho}(K+1)^{(1+\rho)/2}}
 \le\frac{2^{\rho+1}3^{(1-\rho)/2}}{K^{(1+3\rho)/2}},$
where the second inequality uses $K+2\le3K$.  
\end{proof}

\subsection{Proof of Theorem~\ref{thm:pmm}}
When $\alpha_k=1$, we have $\zbf^k=\xbf^{k-1}$, $\tilde{\ybf}^k=\xbf^k$ and $\ybf^k$ such that $v(\ybf^k,f^*)=\min_{1\le s\le k} v(\xbf^s,f^*)$. 
Similar to \eqref{eq:opt-01}, we have 
$f(\xbf^{k})\le\ell_{f}(\xbf^{k},\xbf^{k-1})+\frac{M_{0}(\xbf^{k},\xbf^{k-1})}{1+\rho}\norm{\xbf^{k}-\xbf^{k-1}}^{1+\rho} \le f^{*}+\frac{M_{0}(\xbf^{k},\xbf^{k-1})}{1+\rho}\norm{\xbf^{k}-\xbf^{k-1}}^{1+\rho},$
and $
g_{i}(\xbf^{k})\le\frac{M_{i}(\xbf^{k},\xbf^{k-1})}{1+\rho}\norm{\xbf^{k}-\xbf^{k-1}}^{1+\rho},\ i=1,2,\ldots,m.$
This gives $v(\xbf^{k},f^{*})\le\frac{\hat{M}}{1+\rho}\norm{\xbf^{k}-\xbf^{k-1}}^{1+\rho}.$
Summing up the above result over $k=1,2,\ldots,K$, we have
\begin{equation*}\begin{aligned}
K\min_{1\le k\le K}v(\xbf^{k},f^{*}) & \le\sum_{k=1}^{K}v(\xbf^{k},f^{*})
 \le\sum_{k=1}^{K}\frac{\hat{M}}{1+\rho}\norm{\xbf^{k}-\xbf^{k-1}}^{1+\rho}\\
 & \le\frac{\hat{M}}{1+\rho}\norm{\onebf}_{\frac{2}{1-\rho}}\left(\sum_{k=1}^{K}\norm{\xbf^{k}-\xbf^{k-1}}^{2}\right)^{\frac{\rho+1}{2}}
 \le\frac{\hat{M}_{}}{1+\rho}K^{\frac{1-\rho}{2}}\|\xbf^{*}-\xbf^{0}\|^{\rho+1},
\end{aligned}\end{equation*}
which gives us the desired result. 

\subsection{Proof of Theorem~\ref{thm:restart_rate}}
First, we claim that the inner loop can be terminated at most 
$\blue{K_{s+1}}=\Big\lceil \mathfrak{C}\cdot\theta^{\frac{2(\rho+1-\tilde{\rho})s-2\tilde{\rho}}{\tilde{\rho}(1+3\rho)}}\Big\rceil $
iterations for generating $\blue{\bar{\xbf}^{s+1}}$ that satisfies $\max\{f(\blue{\bar{\xbf}^{s+1}})-f^{*},\norm{[\gbf(\blue{\bar{\xbf}^{s+1}})]_{+}}_{\infty}\}\leq\Delta_{0}\theta^{s+1}$.
We shall prove this result by induction. By definition, we have $\max\{f(\blue{\bar{\xbf}^{0}})-f^{*},\norm{[\gbf(\blue{\bar{\xbf}^{0}})]_{+}}_{\infty}\}=\Delta_{0}$. \blue{Hence, after $K_1$ iterations, we have $\max\bcbra{f(\bar{\xbf}^1)-f^*,\norm{[g(\bar{\xbf}^1)]_+}_{\infty}}\leq \Delta_0 \cdot \theta$ by Theorem~\ref{thm:rate-APM}, which implies that the base case holds.}
Assume at the $s-1$ stage, the output $\blue{\bar{\xbf}^{s}}$ of \apm{} satisfies
$\max\{f(\blue{\bar{\xbf}^{s}})-f^{*},\norm{[\gbf(\blue{\bar{\xbf}^{s}})]_{+}}_{\infty}\}\leq\Delta_{0}\theta^{s}$. 
Then, we consider the upper bound of $\max\{f(\blue{\xbf}^{s+1})-f^{*},\norm{[\gbf(\blue{\xbf}^{s+1})]_+}_{\infty}\}$.
It follows from Theorem~\ref{thm:rate-APM} that 
\begin{equation}
\max\{f(\blue{\bar{\xbf}^{s+1}})-f^{*},\norm{[\gbf(\blue{\bar{\xbf}^{s+1}})]_+}_{\infty}\}\leq\frac{\hat{M}}{1+\rho}\left[{\norm{\proj_{\mathcal{X}^{*}}\cbra{\blue{\bar{\xbf}^{s}}}-\blue{\bar{\xbf}^{s}}}^{\tilde{\rho}}}\right]^{\frac{\rho+1}{\tilde{\rho}}}\frac{2^{\rho+1}3^{(1-\rho)/2}}{\blue{K_{s+1}^{(1+3\rho)/2}}}.\label{eq:des_one_epoch}
\end{equation}
Since $f$ has a H\"{o}lderian growth, we have 
\begin{equation}\label{eq:induction_growth} 
\frac{\mu}{2}\norm{\blue{\bar{\xbf}^{s}}-\proj_{\mathcal{X}^{*}}\cbra{\blue{\bar{\xbf}^{s}}}}^{\tilde{\rho}} \leq f(\blue{\bar{\xbf}^{s}})-f^{*}\leq \max\{f(\blue{\bar{\xbf}^{s}})-f^{*},\norm{[\gbf(\blue{\bar{\xbf}^{s}})]_+}_{\infty}\}\leq\Delta_{0}\cdot\theta^{s}.
\end{equation}
Combining~\eqref{eq:des_one_epoch} and~\eqref{eq:induction_growth}
yields 
\begin{equation}
\max\{f(\blue{\bar{\xbf}^{s+1}})-f^{*},\norm{[\gbf(\blue{\bar{\xbf}^{s+1}})]_+}_{\infty}\}\leq\frac{\hat{M}}{1+\rho}\brbra{\frac{2}{\mu}\Delta_{0}\theta^{s}}^{\frac{\rho+1}{\tilde{\rho}}}\frac{2^{\rho+1}3^{(1-\rho)/2}}{\blue{K_{s+1}^{(1+3\rho)/2}}}.
\end{equation}
\blue{Combining the above result with the definition of $K_{s}$, we obtain $\max\{f(\blue{\bar{\xbf}^{s+1}})-f^{*},\norm{[\gbf(\bar{\xbf}^{s+1})]_+}_{\infty}\}\leq\Delta_{0}\cdot\theta^{s+1}$, completing the induction proof. \blue{Therefore}, to find an $\vep$-optimal solution, we require $\lceil\log_{1/\theta}(\Delta_{0}/\vep)\rceil$ epochs. This implies that the total number of iterations needed by the algorithm to obtain an $\vep$-optimal solution is bounded by $\sum_{s=0}^{\lceil\log_{1/\theta}(\Delta_{0}/\vep)\rceil}K_{s+1} = \sum_{s=0}^{\lceil \log_{1/\theta}(\Delta_{0}/\vep)\rceil}\big\lceil \mathfrak{C}\cdot\theta^{\frac{2(\rho+1-\tilde{\rho})s-2\tilde{\rho}}{\tilde{\rho}(1+3\rho)}}\big\rceil$.}

\blue{For the case $\tilde{\rho}>(1+\rho)$, let $r_{1}=\frac{2(\rho+1-\tilde{\rho})}{\tilde{\rho}(1+3\rho)}<0$ and $r_{2}=-\frac{2}{3\rho+1}$. Then we have the total number of iterations needed by the algorithm to obtain an $\vep$-optimal solution is bounded by:}
\begin{equation}
\begin{aligned} 
\blue{\sum_{s=0}^{\lceil\log_{1/\theta}(\Delta_{0}/\vep)\rceil}\Big\lceil\mathfrak{C}\cdot\theta^{r_{1}s+r_{2}}\Big\rceil}
& \blue{\leq\mathfrak{C}\theta^{r_{2}}\cdot\frac{1-\theta^{r_{1}(2+\log_{1/\theta}\frac{\Delta_{0}}{\vep})}}{1-\theta^{r_{1}}}+2+\log_{1/\theta}\left(\frac{\Delta_{0}}{\vep}\right)}\\
& \blue{=\mathfrak{C}\theta^{r_{2}}\cdot\frac{1-\theta^{2r_{1}}\cdot\left(\frac{\vep}{\Delta_{0}}\right)^{r_{1}}}{1-\theta^{r_{1}}}+2+\log_{1/\theta}\left(\frac{\Delta_{0}}{\vep}\right)}\\
& \blue{=\mathfrak{C}\theta^{r_{2}}\cdot\frac{\theta^{2r_{1}}\cdot\left(\frac{\vep}{\Delta_{0}}\right)^{r_{1}}-1}{\theta^{r_{1}}-1}+2+\log_{1/\theta}\left(\frac{\Delta_{0}}{\vep}\right)}\\
& \blue{\leq\mathfrak{C}\theta^{r_{2}}\cdot(\theta^{r_{1}}-1)^{-1}\cdot\left(\frac{\Delta_{0}}{\theta^{2}\vep}\right)^{-r_{1}}+2+\log_{1/\theta}\left(\frac{\Delta_{0}}{\vep}\right)},
\end{aligned}
\label{eq:upp:rho_leq1}
\end{equation}
\blue{where the last inequality holds since $r_{1}<0$ and $\theta<1$ implies $\theta^{r_{1}}-1>0$.}

For the case $\tilde{\rho}=1+\rho$, we have:
$\blue{\sum_{s=0}^{\lceil \log_{1/\theta}(\Delta_{0}/\vep)\rceil}\Big\lceil \mathfrak{C}\cdot\theta^{-\frac{2}{3\rho+1}}\Big\rceil \leq\Brbra{\mathfrak{C}\theta^{-\frac{2}{3\rho+1}}+1}\cdot\Brbra{\log_{1/\theta}\left(\frac{\Delta_{0}}{\vep}\right)+2}.}$


\section{Missing Proofs of Section~\ref{sec:root-finding}}

\subsection{Proof of Proposition~\ref{prop:value-function}}
The proof of parts 1. and 2. can be found in Section~2.3.4~\citep{nesterov2018lectures}. Part 3 naturally follows from the monotonicity of the directional derivative of a convex function.  \blue{In view of Theorem 3.61 in~\cite{beck2017first}, the 1-Lipschitz continuity of $V(\eta)$ implies $\abs{V^\prime(\eta)}\leq 1$. Since $V^{\prime}(\eta)\leq 0$, then we have $V^{\prime}{(\eta)}\geq -1$.}

Part 4. It is immediate to see $\bar{g}> 0$; otherwise, there would exist an $\xbf\in\bar{\Xcal}$ feasible to \eqref{pb:func-constraint}, which would imply $f^*=\bar{f}$, leading to a contradiction.  Next, we show that for sufficiently small $\eta$, $V(\eta)$ is a linear function. Let $\eta \le \bar{f}-\bar{g}$. It is easy to observe that $
V(\eta)=\min_{\xbf\in\Xcal}\max\Bcbra{f(\xbf)-\eta, \max_{1\le i\le m} g_i(\xbf)}\le \max\{\bar{f}-\eta, \bar{g}\}= \bar{f}-\eta.$
For the reverse direction, we note that $V(\eta)\ge \min_{\xbf\in\Xcal} f(\xbf)-\eta = \bar{f}-\eta$. Thus, we conclude that $V(\eta)=\bar{f}-\eta$.

For part 5,  we have two cases. If $\tilde{\xbf}\in\argmin_{\xbf\in\Xcal} f(\xbf)$, namely, $f(\tilde{\xbf}) = \bar{f}$, by the assumption that $\bar{f}<f^*$, we conclude that $\tilde\xbf$ must be infeasible for \eqref{pb:func-constraint}. Next, we consider $f(\tilde{\xbf})>\bar{f}$.
We observe
\begin{equation}\label{eq:strict-positive}
\begin{aligned}
  v(\tilde\xbf,\eta) & = \min_{\xbf\in\Xcal}\, \max\{f(\xbf)-\eta,g_1(\xbf),\ldots, g_m(\xbf)\}  \\
 & >  \min_{\xbf\in\Xcal}\, \max\{f(\xbf)-\blue{f^*},g_1(\xbf),\ldots, g_m(\xbf)\} = v(\xbf^*,\blue{f^*}) = 0.
\end{aligned}
\end{equation}
In the above, strict inequality holds because, otherwise, $\eta$ is the optimal value. 

We prove the infeasibility of $\tilde{\xbf}$ by contradiction. Suppose that $\tilde{\xbf}$ is feasible, i.e., $\max_{1\le i\le m} g_i(\tilde{\xbf})\le 0$. In view of~\eqref{eq:strict-positive}, we have $f(\tilde{\xbf})>\eta$.
Let $\bar{\xbf}$ be a minimizer of $f(\cdot)$ over the domain $\Xcal$. Clearly, the assumption $\bar{f}<f^*$ implies that $\bar\xbf$ is infeasible, namely, $\max_{1\le i\le m} g_i(\bar\xbf)>0$.
Define $\xbf_\theta = \theta\bar\xbf + (1-\theta)\tilde\xbf$, where $\theta\in [0,1]$. Convexity implies $
f(\xbf_\theta) \le \theta f(\bar{\xbf})+(1-\theta)f(\tilde\xbf) < f(\tilde\xbf), \ \forall \theta \in (0,1],$ 
where the rightmost inequality follows from the assumption $f(\tilde{\xbf})>\bar{f}$. Denote $\delta=f(\tilde\xbf)-\eta-\max_{1\le i\le m} g_i(\tilde\xbf)$. We have $\delta>0$. Since $g_i(\xbf)$ is continuous, for sufficiently small $\theta$, we have $g_i(\xbf_\theta) < g_i(\tilde{\xbf})+\delta$. It follows that \[v(\xbf_\theta,\eta) = \max\{f(\xbf_\theta)-\eta,  \max_{1\le i\le m} g_i(\xbf_\theta)\} < \max\{f(\tilde\xbf)-\eta, \max_{1\le i\le m} g_i(\tilde\xbf)+\delta\} = f(\tilde\xbf)-\eta =v(\tilde\xbf,\eta),\]
which contradicts the optimality of $\tilde\xbf$. Hence, we conclude that $\tilde\xbf$ is infeasible for \eqref{pb:func-constraint}.


\subsection{Proof of Theorem~\ref{thm:convergence-fixedpoint-abstract}}
We show $\eta_{t}\le\blue{f^*}$, $t=0, 1,2,\ldots,$ by induction.
First, we have $\eta_0\le \blue{f^*}$ by our assumption. 
Suppose that $\eta_{s}\le \blue{f^*}$, $s=0,1,\ldots,t-1$. Using the definition of $\eta_{t}$
and the criterion~\eqref{eq:subprob-relative-error}, we have $\eta_{t}=\eta_{t-1}+\beta l_{t-1}\le\eta_{t-1}+\beta V(\eta_{t-1}).$
Moreover, since $0=V(\blue{f^*})\ge V(\eta_{t-1})+V'(\eta_{t-1})(\blue{f^*}-\eta_{t-1})$
and $|V'(\eta)|\le 1$, we have
\begin{equation}
V(\eta_{t-1})\le|V'(\eta_{t-1})|\cdot(\blue{f^*}-\eta_{t-1})\le\blue{f^*}-\eta_{t-1}.\label{eq:v-etak-bound-1}
\end{equation}
Combining these two results, we have
$\eta_{t}\le\brbra{1-\beta }\eta_{t-1}+\beta \blue{f^*}.$
Since $0<\beta <1$, we have $\eta_{t}\le\blue{f^*}$.

Next, we show the linear convergence to optimal level $\blue{f^*}$.
\blue{Since $V$ is a real-valued convex function, the minimum subgradient $\bar{V}^{\prime}$ is well-defined. See Theorem 3.44~\citep{mordukhovich2022convex}.}
Following from the update of $\{\eta_{t}\}$, we have 
\begin{equation}\label{eq:opt-18}
\begin{aligned}
    \blue{f^*}-\eta_{t} & =\blue{f^*}-\eta_{t-1}-\beta l_{t-1} \le\blue{f^*}-\eta_{t-1}-\beta \alpha^{-1}u_{t-1} \le \blue{f^*}-\eta_{t-1}-\beta \alpha^{-1}V(\eta_{t-1}) \\
 &\le\blue{f^*}-\eta_{t-1}-\beta \alpha^{-1}\bsbra{V(\blue{f^*})+\blue{\bar{V}^{\prime}}\cdot(\eta_{t-1}-\blue{f^*})} \\
 & =(1+\beta \alpha^{-1}\blue{\bar{V}^{\prime}})(\blue{f^*}-\eta_{t-1}).
\end{aligned}
\end{equation}
Here, the first inequality is due to \eqref{eq:subprob-relative-error}
and the last inequality follows from the convexity of $V(\cdot)$. 
Applying the convexity of $V(\cdot)$,
\blue{using} \eqref{eq:v-etak-bound-1}, \eqref{eq:subprob-relative-error}, \eqref{eq:opt-18} and $\blue{f^*}-\eta_{0} \le \frac{V(\eta_{0})}{-\blue{\bar{V}^{\prime}}}$, we have \[u_{t}\le\alpha l_{t} \le\alpha V(\eta_{t})
 \le\alpha(\blue{f^*}-\eta_{t})
 \le\alpha\brbra{1+\beta \alpha^{-1}\blue{\bar{V}^{\prime}}}^{t}(\blue{f^*}-\eta_{0})
 \le\exp\brbra{\beta \alpha^{-1}\blue{\bar{V}^{\prime}}\cdot t}\frac{\alpha V(\eta_0)}{-\blue{\bar{V}^{\prime}}}.\]
Setting $\exp\brbra{\beta \alpha^{-1}\blue{\bar{V}^{\prime}}\cdot t}\frac{\alpha V(\eta_0)}{-\blue{\bar{V}^{\prime}}}\le \vep$, we immediately obtain the desired complexity bound.


\subsection{Proof of Theorem~\ref{thm:V-deri}}
For simplicity, let us denote $\tilde{g}_{0}(\xbf,\eta)=f(\xbf)-\eta + \iota_\Xcal(\xbf)$  and $\tilde{g}_{i}(\xbf,\eta)=g_{i}(\xbf)+ \iota_\Xcal(\xbf)$ ($1\le i \le m$), where $\iota_\Xcal$ is the indicator function on $\Xcal$. See the definition in Section~\ref{subsec:prelim}.  Denote $I_m=\{i: \tilde{g}_i(\xbf,\eta)=v(\xbf,\eta), 0\le i\le m\}$.
In view of Theorem~3.59~\citep{mordukhovich2022convex}, we obtain
 \begin{equation}\label{eq:subdiff-v}
     \partial v(\xbf,\eta)\supseteq \textrm{co}\big\{\cup_{i\in I_m}\partial\tilde{g}_{i}(\xbf,\eta)\big\},
 \end{equation}
 where $\textrm{co}$ denotes the convex hull and 
$\textrm{co}\left(\cup_{i\in I_{m}}\partial\tilde{g}_{i}(\xbf,\eta)\right)=\left\{ \tsum_{i\in I_{m}}\lambda_i \partial \tilde{g}_i(\xbf,\eta):\, \tsum_{i\in I_{m}} \lambda_i =1, \lambda_i\ge 0\right\}.$ 
Note that the equality in~\eqref{eq:subdiff-v} holds when $\xbf \in \textrm{int}(\Xcal)$.
We compute the subdifferentials $\partial \tilde{g}_{0}(\xbf,\eta)=[\partial f(\xbf) + N_\Xcal(\xbf)]\times \{-1\}$, and $\partial\tilde{g}_{i}(\xbf,\eta)=[\partial g_i(\xbf) + N_\Xcal(\xbf)]\times \{0\}$. \blue{The exchange of subdifferentials and summations is guaranteed by the Moreau-Rockafellar sum rule~\citep[Thm 2.26]{royset2021optimization}.}
Define $\Lambda(\xbf, \eta)=\big\{\lambf=[\lambda_0,\lambda_1,\ldots,\lambda_m]^\top\in\Sbb^{m}: \lambda_0=0 \text{ if } f(\xbf)-\eta<v(\xbf,\eta), \text{ and } \lambda_i =0 \text{ if } {g_i}(\xbf,\eta) < v(\xbf,\eta),1\le i\le m\big\}$. Here, $\Sbb^m$ is the standard simplex of dimension $m$.
Then the term $\textrm{co}\left(\cup_{i\in I_{m}}\partial\tilde{g}_{i}(\xbf,\eta)\right)$ can be expressed as:
\begin{equation}\label{eq:convex-hull}
    \textrm{co}\left(\cup_{i\in I_{m}}\partial\tilde{g}_{i}(\xbf,\eta)\right)
    = \Big\{ \Brbra{\lambda_{0}\partial f(\xbf)+\sum_{j=1}^m\lambda_{j}\partial g_{j}(\xbf)+\mcal N_\Xcal(\xbf)}\times\left(-\lambda_{0}\right):\lambf\in\Lambda(\xbf,\blue{\eta})\Big\}.
\end{equation}

 Based on \citet[Theorem 10.13]{rockafellar2009variational}, we 
express the subdifferential  $V(\eta)$ by
\begin{equation}\label{eq:subdiff-Veta}
\partial V(\eta)=\cup_{\xbf\in S(\eta)}M(\xbf,\eta),  \ \text{where } M(\xbf,\eta)=\bcbra{y:(0,y)\in\partial v(\xbf,\eta)},
\end{equation} 
where $S(\eta)=\argmin_{\xbf\in\Xcal }v(\xbf,\eta)$ is the solution set. 
In view of \eqref{eq:subdiff-v}, \eqref{eq:convex-hull} and \eqref{eq:subdiff-Veta}, we have 
\[
\partial V(\eta) \supseteq \bigcup_{\xbf\in S(\eta)} \big\{ -\lambda_0: \ 0\in \lambda_{0} \partial f(\xbf)+\sum_{j=1}^m\lambda_{j} \partial g_{j}(\xbf)+\mcal N_\Xcal(\xbf),\,\lambf\in\Lambda(\xbf, \eta) \big\}.
\]
We fix $\eta=f^*$ and $\xbf^{*}\in S(f^*)$. 
 Note that the condition
\begin{equation}\label{eq:stationary}
\begin{aligned}
&0\in \lambda_{0}\partial f(\xbf^{*})+\sum_{j = 1}^{m}\lambda_{j}\partial g_{j}(\xbf^{*})+\mcal N_\Xcal(\xbf^*), \, \lambf\in\Lambda(\xbf^*, \blue{f^*})
\end{aligned}
\end{equation}
characterizes the Fritz-John condition \blue{of problem~\eqref{pb:func-constraint}} at optimality. 

\blue{Under a standard constraint qualification (e.g., Slater's condition), the KKT condition of~\eqref{pb:func-constraint} holds and we have $\lambda_0>0$. We transition from Fritz-John to the KKT condition by dividing \eqref{eq:stationary} by $\lambda_0$. We define the Lagrange multiplier vector $\ybf^*\in \Rbb_+^m$ as $y^*_i = \frac{\lambda_i}{\lambda_0}$, $1\le i\le m$. From the simplex constraint $\tsum_{i=0}^m \lambda_i=1$, it follows that $\lambda_0(1+\tsum_{i=1}^m y_i^*)=1$. This implies $\lambda_{0}={1}/(1+\norm{\ybf^*}_{1})$, which completes our proof.}

\subsection{Proof of Theorem~\ref{thm:convergence-secant-1}}
The monotonicity of $\{\eta_t\}$ can be easily derived based on the non-negativity of stepsizes and $l_t$.
    For brevity, we denote $w_t\coloneqq\frac{u_{t-1}-l_t}{\eta_{t-1}-\eta_t}$, for $t=1,2,3,\ldots$. 
    We show that the monotonicity of $\{\eta_t\}$ and $\eta_t\le \blue{f^*}$ by induction. 
    Suppose $\eta_s\le \blue{f^*}$ for $s=0, 1,2,\ldots, t-1$, and $\{\eta_s\}_{0\le s \le t-1}$ is monotonically increasing.
    Using the convexity of $V(\eta)$, we have \[w_s\le \frac{V(\eta_{s-1})-V(\eta_s)}{\eta_{s-1}-\eta_s}\le V'(\eta_s)<0, \ s=1,2,\ldots, t-1.\]
    We consider two cases.  First, suppose that $w_{t-1}\ge -1$.  Then,
    \begin{equation}
        \begin{aligned}
            \eta_t  &= \eta_{t-1} - \beta \frac{l_{t-1}}{w_{t-1}} \le \eta_{t-1} - \beta \frac{V(\eta_{t-1})}{V'(\eta_{t-1})}
             \le \eta_{t-1} -\beta \frac{V'(\eta_{t-1})(\eta_{t-1}-\blue{f^*})}{V'(\eta_{t-1})}  =(1-\beta)\eta_{t-1} + \beta \blue{f^*} \le \blue{f^*}.
        \end{aligned}
    \end{equation}
    Second, suppose $w_{t-1}<-1$, then Algorithm~\ref{alg:inexact-secant} reduces to fixed point iteration, and we have $\eta_t\le \blue{f^*}$ by the proof of Theorem~\ref{thm:convergence-fixedpoint-abstract}.  

    Next, we bound the total number of iterations. 
    Let $\Delta_t = \blue{f^*}-\eta_t$. Then the update in line~\ref{tag:eta-update-secant} of Algorithm~\ref{alg:inexact-secant} reads \[\Delta_t = \Delta_{t-1} + \beta \min \bcbra{-l_{t-1}, \frac{l_{t-1}}{u_{t-2}-l_{t-1}}(\Delta_{t-1}-\Delta_{t-2}) }.\]
    Using the first term in the $\min$ function, we have
    $\Delta_t \le \Delta_{t-1}  -\beta  l_{t-1} $. Similar to the argument showing \eqref{eq:opt-18}, we obtain $\blue{f^*}-\eta_{t}\le\sigma(\blue{f^*}-\eta_{t-1})$ and the linear rate: $
    \blue{f^*}-\eta_t \le \sigma^{t-1} (\blue{f^*}-\eta_1) \le \sigma^{t} (\blue{f^*}-\eta_0) , \forall t=1,2,3,\ldots$
    Moreover, by applying the second part of the $\min$ function, we have \[\Delta_t \le  \Delta_{t-1} + \beta \frac{l_{t-1}}{u_{t-2}-l_{t-1}}(\Delta_{t-1}-\Delta_{t-2}) \le \Delta_{t-1} + \beta \frac{l_{t-1}}{\alpha l_{t-2}-l_{t-1}}(\Delta_{t-1}-\Delta_{t-2}),\]
    where the last inequality follows from \eqref{eq:subprob-relative-error}.
    Let us denote $\lambda_j=\frac{l_j}{\alpha l_{j-1}}$, then we have $\beta \frac{\lambda_{\blue{j}-1}}{1-\lambda_{\blue{j}-1}} \le \frac{\Delta_{\blue{j}}-\Delta_{\blue{j}-1}}{\Delta_{\blue{j}-1}-\Delta_{\blue{j}-2}}.$
    Multiply the above relation for $\blue{j}=2,3,\ldots,t+1$, we have
     \begin{equation}\label{eq:mid-02}
    \beta^t \frac{\prod_{j=1}^t\lambda_{j}}{\prod_{j=1}^t(1-\lambda_{j})} \le \frac{\Delta_{t+1}-\Delta_{t}}{\Delta_{1}-\Delta_0}.
    \end{equation}
    Notice that the relation \blue{$l_{j}<u_{j}\leq u_{j-1}\leq\alpha l_{j-1}$ implies $\lambda_j <1$}. Moreover, we have $\prod_{j=1}^t\lambda_{j}\cdot \prod_{j=1}^t(1-\lambda_{j}) \le (\frac{1}{4})^t$, which uses the relation $2 [\lambda_j (1-\lambda_j)]^{1/2} \le \lambda_j +(1-\lambda_j) = 1$ from the arithmetic mean inequality. 
    Multiplying this relation with \eqref{eq:mid-02}, we have 
    $\beta^t\left(\frac{1}{\alpha^t} \frac{l_t}{l_0}\right)^2 = \beta^t \Brbra{\prod_{j=1}^t\lambda_{j}}^2 \le (\frac{1}{4})^t \frac{\Delta_{t+1}-\Delta_{t}}{\Delta_{1}-\Delta_0}\le (\frac{1}{4})^t \frac{\blue{f^*}-\eta_0}{\eta_1-\eta_0},
    $
    which implies $l_t \le l_0 \sqrt{\frac{\blue{f^*}-\eta_0}{\eta_1-\eta_0}}\brbra{\frac{\alpha}{2\sqrt{\beta}}}^t$. Setting $T=\log_{2\sqrt{\beta}/\alpha}\big(\frac{\alpha l_0}{\varepsilon} \sqrt{\frac{\blue{f^*}-\eta_0}{\eta_1-\eta_0}}\big)$, we have $u_T\le \alpha l_T\le \vep$.

\section{Convergence analysis in Section~\ref{sec:apl-based}\label{sec:convergence_ana_sec_4}}

\subsection{More details about the APL method}\label{subsec:gap-reduce}
The following result summarizes some important convergence properties of the gap reduction.
\begin{proposition}\label{prop:gap-reduction}
Problem~\eqref{eq:update-xk-2} is always feasible unless it terminates at line~\ref{tag:gap-reduce-break-1}. Whenever the algorithm terminates, the latest estimated bounds $\tilde{v}_k^{\rm{U}}$ and $\tilde{v}_k^{\rm{L}}$ satisfy:
$\tilde{v}_k^{\rm{U}} - \tilde{v}_k^{\rm{L}} \le \brbra{\blue{\frac{1+\theta}{2}}}(\tilde{u}-\tilde{l}).$
Suppose the algorithm computes $\xbf^1,\ldots, \xbf^K$, then we have
\begin{equation}\label{eq:opt-05}
\Gamma_{K}[v(\ybf^{K},\eta)-\lambda] \le\Gamma_{1}(1-\alpha_{1})[v(\ybf^{0},\eta)-\lambda]+\frac{\bar{M}}{1+\rho}\norm{\cbf_{K}}_{\frac{2}{1-\rho}}\left(\norm{\xbf^{K}-\xbf^{0}}^{2}\right)^{\frac{\rho+1}{2}},
\end{equation}
where $\Gamma_K$, $\cbf_K$ is defined in Theorem~\ref{thm:rate-APM} and $\bar{M}=\max_{0\le i\le m}\sup_{\ybf,\zbf\in\Xcal}M_i(\ybf,\zbf)$.
Moreover, if we set $\alpha_k = 2/(k+1)$, then the iteration number of {APL gap reduction} is bounded by
$\hat{K} = \bigg\lceil\Big(\frac{2^{\rho+1}3^{(1-\rho)/2}\bar{M}D_\Xcal^{(\rho+1)/2}}{(1+\rho)\theta (\tilde{u}-\tilde{l})}\Big)^{{2}/{(1+3\rho)}}\bigg\rceil.$
\end{proposition}
\begin{proof}
For $k>1$, since $\xbf^k\in \bar{\Xcal}_{k-1}\subseteq \bar{\Xcal}_{k-1}^{\rm U}$, by definition of $\mcal X_{k-1}^U$, we have 
$\inner{\xbf^{k-1}-\xbf^0}{\xbf^k-\xbf^{k-1}} \ge 0$,
which implies $\norm{\xbf^k-\xbf^0}^2 -\norm{\xbf^{k-1}-\xbf^0}^2  - \norm{\xbf^k-\xbf^{k-1}}^2 \ge 0. $
Note that this property naturally holds when $k=1$.
Summing up the above relation for $k=K, K-1,\ldots, 1$, we have
\begin{equation}\label{eq:opt-07}
\sum_{k=1}^K \norm{\xbf^k-\xbf^{k-1}}^2 \le \norm{\xbf^K-\xbf^0}^2.
\end{equation}

Suppose that subproblem~\eqref{eq:update-xk-2} is infeasible at the $K$-th iteration, then we have $h_K>\lambda$, and hence $v_K^{\rm L}\blue{\geq} \lambda$. Then the break condition in line~\ref{tag:gap-reduce-break-1} is met. 
     Otherwise, suppose subproblem~\eqref{eq:update-xk-2} is feasible for $k=1,2,\ldots, K$. Similar to the proof of Theorem~\ref{thm:rate-APM}, we can show $
\Gamma_{K}[v(\ybf^{k},\eta)-\lambda] \le\Gamma_{1}(1-\alpha_{1})[v(\ybf^{0},\eta)-\lambda]+\frac{\hat{M}}{1+\rho}\norm{\cbf_{K}}_{\frac{2}{1-\rho}}\left(\sum_{k=1}^{K}\norm{\xbf^{k}-\xbf^{k-1}}^{2}\right)^{(\rho+1)/2}.$ 
Combining this result with \eqref{eq:opt-07} gives \eqref{eq:opt-05}.

Setting $\alpha_k=2/(k+1)$ and using the fact that $D_{\Xcal}\ge \norm{\xbf^K-\xbf^0}^2$, we have $
v(\ybf^K,\eta)-\lambda \le \frac{\bar{M}}{1+\rho}\left[D_{\Xcal}\right]^{\frac{\rho+1}{2}}\frac{2^{\rho+1}3^{(1-\rho)/2}}{K^{(1+3\rho)/2}}.$
It is easy to check that when the iteration number reaches $\hat K$, the termination criterion of line~\eqref{tag:gap-reduce-break-2} must be satisfied.

\end{proof}

\subsection{\blue{Proof of Theorem~\ref{thm:complexity-APL}}}
The proof is adapted from \citet{lan2015bundle} to accommodate our termination criteria.

Let $\delta_s=\bar{u}_s-\bar{l}_s$. We have the property that $\delta_{s+1}\le \gamma \delta_s$.
Applying Proposition~\ref{prop:gap-reduction}, we have that at the $s$-th stage, the gap reduction will terminate in at most 
$$\hat{K}_s = \bigg\lceil\left(\frac{2^{\rho+1}3^{(1-\rho)/2}\bar{M}D_\Xcal^{(\rho+1)/2}}{(1+\rho)\theta (\tilde{u}_{s-1}-\tilde{l}_{s-1})}\right)^{{2}/{(1+3\rho)}}\bigg\rceil=\mcal O\Brbra{\Big\lceil\brbra{\bar{u}_{s-1}-\bar{l}_{s-1}}^{-\frac{2}{1+3\rho}}\Big\rceil} = O\Brbra{\Big\lceil\delta_{s-1}^{-\frac{2}{1+3\rho}}\Big\rceil}$$ 
iterations.
To bound the total number of calls to the gap reduction, it suffices to establish a bound on $\sum_{s=1}^S \Big\lceil \delta_{s-1}^{-\frac{2}{1+3\rho}} \Big\rceil$. We proceed by considering the following two cases. 

\textbf{Case 1:} Suppose $V(\eta)\geq \frac{1}{2}\vep$. 
{Note that for any $s>0$, we have $\delta_s=\bar{u}_s-\bar{l}_s\leq \gamma^{s}(\bar{u}_0 - \bar{l}_0) $ and $ \frac{\alpha-1}{\alpha}V(\eta)\le \frac{\alpha-1}{\alpha} \bar{u}_s.$
    Letting 
    $ S_1= \max\left\{0, \left\lceil\log_{1/\gamma} \frac{\alpha(\bar{u}_0-\bar{l}_0)}{(\alpha-1)V(\eta)}\right\rceil\right\}$, we can see that $\gamma^{S_1}\brbra{\bar{u}_0 - \bar{l}_0}\leq \frac{\alpha-1}{\alpha}V(\eta)$, which implies $\bar{u}_{S_1} - \bar{l}_{S_1}\leq \frac{\alpha - 1}{\alpha} \bar{u}_{S_1}$ and satisfies the first condition of the while loop.  Hence, it takes at most $S_1$ iterations to terminate the while loop.} Without loss of generality, we assume $\delta_0>\frac{\alpha-1}{\alpha} \bar{u}_0$ and $\bar{u}_0>\vep$. 
    Observe that
 $\frac{\alpha(\bar{u}_0-\bar{l}_0)}{(\alpha-1)V(\eta)}
\;\ge\;
\frac{\alpha\,\delta_0}{(\alpha-1)\bar{u}_0}
\;>\;1.$ This implies ${S_1>0}$.
    Consequently, we have
    \begin{equation}\label{delta-sum-1}
        \sum_{s=1}^{{S_1}} \Big\lceil \delta_{s-1}^{-\frac{2}{1+3\rho}} \Big\rceil \le {S_1} + \sum_{s=1}^{{S_1}}\Big( \delta_{s-1}^{-\frac{2}{1+3\rho}} \Big) 
         \le {S_1} + \sum_{s=1}^{{S_1}} \Big( \frac{\alpha \gamma^{(S_1-s)}}{(\alpha-1)V(\eta)}\Big)^{\frac{2}{1+3\rho}}
         \le {S_1} + \frac{1}{1-\gamma^{2/(1+3\rho)}}\Big(\frac{\alpha}{(\alpha-1)V(\eta)}\Big)^{\frac{2}{1+3\rho}}.
    \end{equation}
    Note that the bound trivially holds when $S_1=0$.

\textbf{Case 2:} Consider when $V(\eta)<\frac{1}{2}\vep$. 
{Letting
$S_2=\max\left\{0,\left\lceil\log_{1/\gamma}\frac{2(\bar{u}_0 - \bar{l}_0)}{\vep}\right\rceil\right\}$,  we can see that $\bar{u}_{S_2} - \frac{1}{2}\vep \leq \bar{u}_{S_2} - V(\eta)\leq \bar{u}_{S_2} - \bar{l}_{S_2}\leq \gamma^{S_2} (\bar{u}_0 - \bar{l}_0)$ holds. Hence, it takes at most $S_2$ iterations to exit the while loop for the second condition $\bar{u}_{S_2} < \vep$ since $\gamma^{S_2} (\bar{u}_0 - \bar{l}_0)<\frac{1}{2}\vep$. 
    Without loss of generality, we assume 
    $\bar{u}_0>\vep$. Consequently, $\frac{2(\bar{u}_0-\bar{l}_0)}{\vep}\ge \frac{2(\bar{u}_0-\vep/2)}{\vep}>1$. } This implies  ${S_2>0}$.
    According to the while condition of the algorithm, we have $\delta_{S_2-1}>\frac{\alpha-1}{\alpha}\bar{u}_{S_2-1} > \frac{\alpha-1}{\alpha}\vep$. It follows that 
    $\frac{\alpha-1}{\alpha}\vep \le \delta_{S_2-1}\le \ldots \le  \gamma^{S_2-1-s}\delta_s.$
    Consequently, we have
    \begin{equation}\label{delta-sum-2}
        \sum_{s=1}^{{S_2}} \Big\lceil \delta_{s-1}^{-\frac{2}{1+3\rho}} \Big\rceil 
        \le {S_2} + \sum_{s=1}^{{S_2}}\Big( \delta_{s-1}^{-\frac{2}{1+3\rho}} \Big) 
        \le {S_2} + \sum_{s=1}^{{S_2}} \Big( \frac{\alpha \gamma^{(S_2-s)}}{(\alpha-1)\vep}\Big)^{\frac{2}{1+3\rho}} 
        \le {S_2} + \frac{1}{1-\gamma^{2/(1+3\rho)}}\Big(\frac{\alpha}{(\alpha-1)\vep}\Big)^{\frac{2}{1+3\rho}}.
    \end{equation}
Next, we combine \eqref{delta-sum-1} and \eqref{delta-sum-2} to provide a unified bound. Let $\bar{S}_2=\max\left\{0,\left\lceil\log_{1/\gamma}\frac{2\alpha(\bar{u}_0 - \bar{l}_0)}{(\alpha-1)\vep}\right\rceil\right\}$, we have $\bar{S}_2\ge S_2$.
If $V(\eta)\ge \frac{1}{2}\vep$, since $\min\{\frac{1}{V(\eta)}, \frac{2}{\vep}\}=\frac{1}{V(\eta)}$, we have $S_1\le S$. Otherwise, $\min\{\frac{1}{V(\eta)}, \frac{2}{\vep}\}=\frac{2}{\vep}$ and $\bar{S}_2\le S$. 
Using a similar analysis, we can provide a unified bound on the second terms of \eqref{delta-sum-1} and \eqref{delta-sum-2}. Thus, we have
$
\sum_{s=1}^{{S_2}} \Big\lceil \delta_{s-1}^{-\frac{2}{1+3\rho}} \Big\rceil \le S + \frac{1}{1-\gamma^{2/(1+3\rho)}}\Big(\frac{2\alpha}{(\alpha-1)}\min\left\{\frac{1}{V(\eta)}, \frac{1}{\vep}\right\}\Big)^{\frac{2}{1+3\rho}}.
$

\subsection{Proof of Theorem~\ref{thm:fixed_point_secant_all_complexity}\label{sec:EC3}}


We now prove each of the three claims in sequence.

Part 1. \blue{For Algorithm~\ref{alg:prac-ifp}, the non-negativity of $\tilde{l}_t$ is straightforward to see.} 
Due to convexity of $V(\cdot)$, we have $V(\eta_t)\geq V(\eta_{t-1})+V^{\prime}(\eta_{t-1})(\eta_t - \eta_{t-1})$. Since $|V'(\eta)|\le 1$, we have 
\[
V(\eta_{t})\ge V(\eta_{t-1})-|\eta_{t-1}-\eta_{t}|\ge l_{t-1}-\beta l_{t-1}=(1-\beta)l_{t-1}.\]
Furthermore, using the monotonicity property from Proposition~\ref{prop:value-function}, we have 
\[V(\eta_t)
        \geq V(\eta_{t-1})+\frac{V(\eta_{t-1})-V(\eta_{t-2})}{\eta_{t-1}-\eta_{t-2}}\beta l_{t-1}
        \geq l_{t-1}+\frac{l_{t-1}-u_{t-2}}{\beta l_{t-2}}\beta l_{t-1} = (1+\frac{l_{t-1}-u_{t-2}}{l_{t-2}})l_{t-1}, \ \text{for } t\geq 2.\]
\blue{When $t=1$, we have $(1+\frac{l_{t-1}-u_{t-2}}{l_{t-2}})=0$.}
Combining these lower bounds, we see that $\tilde{l}_{t}$ defined in Algorithm~\ref{alg:prac-ifp} is a valid lower bound on $V(\eta_{t})$. 

\blue{For Algorithm~\ref{alg:secant-2}, the non-negativity of $\tilde{l}_t$ is easy to observe.}
\blue{Due to the convexity of $V(\cdot)$, we have $
 V(\eta_t) \ge V(\eta_{t-1})+V'(\eta_{t-1})(\eta_t-\eta_{t-1}).$
 For $t=1$, since $V'(\eta)\ge -1$, it is easy to see
 $ V(\eta_t) \ge V(\eta_{t-1})-\beta l_{t-1}\ge (1-\beta)l_{t-1}$.
For $t>1$, following the update rule and convexity, we have}
\begin{equation}
\begin{aligned}
V(\eta_t) 
& = V(\eta_{t-1})- \beta V'(\eta_{t-1})\frac{\eta_{t-2}-\eta_{t-1}}{u_{t-2}-l_{t-1}}l_{t-1}\\
& \ge V(\eta_{t-1})- \beta V'(\eta_{t-1})\frac{\eta_{t-2}-\eta_{t-1}}{V(\eta_{t-2})-V(\eta_{t-1})}l_{t-1}\\
& \ge \left[1-\beta V'(\eta_{t-1})\frac{\eta_{t-2}-\eta_{t-1}}{V(\eta_{t-2})-V(\eta_{t-1})}\right] l_{t-1}\geq (1-\beta)l_{t-1},
\end{aligned}
\end{equation}
\blue{where the last inequality uses Proposition~\ref{prop:value-function}. Hence, $\tilde{l}_t$ is a valid nonnegative lower bound on $V(\eta_t)$.}

Part 2. 
\blue{For both algorithms, we have $V(\eta_T)\geq \tilde{l}_T > 0 = V(f^*)$. Since $V(\eta)$ is non-increasing, this implies $\eta_T < f^*$. Hence, taking $\eta_1=\eta_T,\eta_2=\blue{f^*}$ in Proposition~\ref{prop:value-function}, we have $\blue{f^*}-\eta_T\leq\frac{V(\eta_T)-V(\blue{f^*})}{-\blue{\bar{V}^{\prime}}}\leq\frac{u_T}{-\blue{\bar{V}^{\prime}}}\le\frac{\vep}{-\blue{\bar{V}^{\prime}}}$. When the algorithm terminates, i.e., $v(\xbf^T;\eta_{T})\leq u_{T}\le\vep$, we have $\norm{[\gbf(\xbf^{T})]_{+}}_{\infty}\le\vep$ and $f(\xbf^{T})-\eta_{T}\le\vep$.
It follows that 
$f(\xbf^{T})-f^* =f(\xbf^{T})-\blue{f^*}\le f(\xbf^{T})-\eta_{T}\le\vep.$}

Part 3. {For both Algorithm~\ref{alg:prac-ifp} and Algorithm~\ref{alg:secant-2}, by Theorem~\ref{thm:complexity-APL}, at the $t$-th iteration the number of gap reduction iterations in APL is at most $K_{t}=\mcal O\Brbra{\log\brbra{\min\{\frac{1}{V(\eta_{t})},\frac{1}{\vep}\}}+\brbra{\min\{\frac{1}{V(\eta_{t})},\frac{1}{\vep}\}}^{\frac{2}{1+3\rho}}}.$ For $t=1,\ldots,T-1$, APL terminates at the condition $u_t\leq \alpha l_t$. Hence, it follows from Theorem~\ref{thm:complexity-APL} that $K_t$ reaches the first bound during the minimization of the two, i.e.,  $K_{t}=\mcal O\Brbra{\log(\frac{1}{V(\eta_{t})})+\brbra{\frac{1}{V(\eta_{t})}}^{\frac{2}{1+3\rho}}}$. Recall that we have $\blue{f^*}-\eta^{T}\leq \sigma^{T-t}(\blue{f^*}-\eta^{t})$ and $\frac{1}{\blue{f^*}-\eta^{T-1}}\leq \frac{-V'(\eta^{T-1})}{V(\eta^{T-1})}\leq \frac{-V'(\eta^{T-1})}{\vep}$ from 
 Proposition~\ref{prop:value-function}.
Due to the convexity of $V(\cdot)$, for any $s\in\{0, 1,2,\ldots, T-1\}$,  we have }
\begin{equation}\label{eq:mid-04}
    \begin{aligned}
    \sum_{t=1}^{T-1}\min\{\frac{1}{V(\eta_{t})},\frac{1}{\vep}\}^{\frac{2}{1+3\rho}}&\le\sum_{t=1}^{T-s-1}\left(\frac{1}{V(\eta_{t})}\right)^{\frac{2}{1+3\rho}}+\sum_{t=T-s}^{T-1}\left(\frac{1}{\vep}\right)^{\frac{2}{1+3\rho}}\\
    &\leq\sum_{t=1}^{T-s-1}\left(\frac{1}{-\blue{\bar{V}^{\prime}} \cdot (\blue{f^*}-\eta_{t})}\right)^{\frac{2}{1+3\rho}}+s\left(\frac{1}{\vep}\right)^{\frac{2}{1+3\rho}}\\
    &\leq\sum_{t=1}^{T-s-1}\left(\frac{\sigma^{T-t-1}}{-\blue{\bar{V}^{\prime}}\cdot (\blue{f^*}-\eta^{T-1})}\right)^{\frac{2}{1+3\rho}}+s\left(\frac{1}{\vep}\right)^{\frac{2}{1+3\rho}},\\
\end{aligned}
\end{equation}
where the second inequality uses $V(\eta_{t})\ge0+\blue{\bar{V}^{\prime}}\cdot (\eta_{t}-\blue{f^*})$, and the last one uses the contraction property:
$\blue{f^*}-\eta^{T}\leq \sigma^{T-t}(\blue{f^*}-\eta^{t})$.
Since $V(\eta_{T-1})-0\le V'(\eta_{T-1})(\eta_{T-1}-\blue{f^*})$ and 
$\alpha V(\eta_{T-1})\geq \alpha l_{T-1}\geq u_{T-1}\geq \vep$, it follows that
\begin{equation}\label{eq:upper_bound_sum}
    \begin{aligned}
    \sum_{t=1}^{T-1}\min\{\frac{1}{V(\eta_{t})},\frac{1}{\vep}\}^{\frac{2}{1+3\rho}}
    &\leq\sum_{t=1}^{T-s-1}\left(\frac{-V^{\prime}(\eta^{T-1})\sigma^{T-t-1}}{-\blue{\bar{V}^{\prime}}\cdot V(\eta_{T-1})}\right)^{\frac{2}{1+3\rho}}+s\left(\frac{1}{\vep}\right)^{\frac{2}{1+3\rho}}\\
    &\leq\sum_{t=1}^{T-s-1}\left(\frac{\alpha V^{\prime}(\eta^{0})\sigma^{T-t-1}}{\blue{\bar{V}^{\prime}}\cdot \vep}\right)^{\frac{2}{1+3\rho}}+s\left(\frac{1}{\vep}\right)^{\frac{2}{1+3\rho}}\\
    &\le\left(\frac{\alpha V'(\eta_{0})}{\blue{\bar{V}^{\prime}}\cdot \vep}\right)^{\frac{2}{1+3\rho}}\sum_{t=1}^{T-s-1}\brbra{\sigma^{\frac{2}{1+3\rho}}}^{T-t-1}+s\left(\frac{1}{\vep}\right)^{\frac{2}{1+3\rho}}\\
    &\le\left(\frac{\alpha V'(\eta_{0})}{\blue{\bar{V}^{\prime}}\cdot \vep}\right)^{\frac{2}{1+3\rho}}\frac{\sigma^{\frac{2}{1+3\rho}s}}{1-\sigma^{\frac{2}{1+3\rho}}}+s\left(\frac{1}{\vep}\right)^{\frac{2}{1+3\rho}}.
\end{aligned}
\end{equation}
Note that the above relation also holds for any $s\ge T$.
We want to set $s$ to obtain a sharper upper bound. 
Let $s^{*}:
=\frac{2/(1+3\rho)\log(-\blue{\bar{V}^{\prime}})+\log\left(1-\sigma^{{2}/{(1+3\rho)}}\right)}{\log\sigma^{{2}/{(1+3\rho)}}}$. We immediately have $\ensuremath{\frac{1}{|\blue{\bar{V}^{\prime}}|^{2/(1+3\rho)}}\frac{\sigma^{{2s^{*}}/{(1+3\rho)}}}{1-\sigma^{{2}/{(1+3\rho)}}}\le1}$. Moreover, 
using $\log(1-x)\le-x$ and $(1-x)^{a}\le1-ax$ for $x>0$ and $a\in(0,1)$, and then applying Theorem~\ref{thm:V-deri}, we have 
\begin{equation}
    \begin{aligned}
    s^{*} &\le\frac{\frac{2}{1+3\rho}\log(-\blue{\bar{V}^{\prime}})+\log\left(-\frac{2}{1+3\rho}\frac{\beta}{\alpha}\blue{\bar{V}^{\prime}}\right)}{\frac{2}{1+3\rho}\frac{\beta}{\alpha}\blue{\bar{V}^{\prime}}}\\
    &=\Bsbra{\frac{(3+3\rho)\alpha}{2\beta}\log\left(1+\norm{\ybf^{*}}\right)+\frac{(1+3\rho)\alpha}{2\beta}\log\left(\frac{(1+3\rho)\alpha}{2\beta}\right)}(1+\norm{\ybf^{*}})=:T_{\ybf^{*}}.
    \end{aligned}
\end{equation}

\blue{Now, we are ready to give the specific complexity of Algorithm~\ref{alg:prac-ifp} and Algorithm~\ref{alg:secant-2}.}

\blue{For Algorithm~\ref{alg:prac-ifp}, }
we note that the iteration number $T$ is bounded by $T\le  T_{\vep}^{{\rm {FP}}}=\Ocal\Brbra{{\frac{\alpha}{-\beta \blue{\bar{V}^{\prime}}}\log\brbra{\frac{\alpha V(\eta_{0})}{-\blue{\bar{V}^{\prime}}\vep}}}}$. 
If $T_{\ybf^{*}}\le T_{\vep}^{{\rm {FP}}}-1$, following \eqref{eq:upper_bound_sum} and setting $s=s^{*}$, we have 
\begin{equation*}
    \begin{aligned}
    & \sum_{t=1}^{T-1}\min\bcbra{\frac{1}{V(\eta_{t})},\frac{1}{\vep}}^{\frac{2}{1+3\rho}}\le\left(\frac{\alpha|V'(\eta_{0})|}{\vep}\right)^{\frac{2}{1+3\rho}}+T_{\ybf^{*}}\left(\frac{1}{\vep}\right)^{\frac{2}{1+3\rho}} \\
    & \quad =\Ocal\Brbra{\left((\norm{\ybf^{*}}_{1}+1)\log(\norm{\ybf^{*}}_{1}+1)\right)\left(\frac{1}{\vep}\right)^{\frac{2}{1+3\rho}}}.
    \end{aligned}
\end{equation*}
If $T_{\ybf^{*}}>T_{\vep}^{{\rm {FP}}}-1$, then $T_\vep^{\rm{FP}}\le \Ocal\brbra{(\norm{\ybf^{*}}_{1}+1)\log\brbra{\norm{\ybf^{*}}_{1}+1}}$.
By setting $s=T-1$, we have 
\begin{equation*}
    \begin{aligned}
    & \sum_{t=1}^{T-1}\min\bcbra{\frac{1}{V(\eta_{t})},\frac{1}{\vep}}^{\frac{2}{1+3\rho}}
    \le (T-1)\left(\frac{1}{\vep}\right)^{\frac{2}{1+3\rho}}
    < T_\vep^{\rm{FP}} \big(\frac{1}{\vep}\big)^{\frac{2}{1+3\rho}} \\
    & \quad =\Ocal\Brbra{\left((\norm{\ybf^{*}}_{1}+1)\log\brbra{\norm{\ybf^{*}}_{1}+1}\right)\left(\frac{1}{\vep}\right)^{\frac{2}{1+3\rho}}}.
    \end{aligned}
    \end{equation*}
Furthermore, since $V(\eta_t) \geq \vep /\alpha$, 
and there exists a $\ybf^*$ such that $\blue{\bar{V}^{\prime}}=-\frac{1}{1+\norm{\ybf^*}}$, we have $\mcal O\Brbra{\sum_{t=1}^{T-1}\log(\frac{1}{V(\eta_{t})})}=\mcal O{\Brbra{\frac{\alpha(1+\norm{\ybf^{*}})}{\beta}\log\brbra{\frac{\alpha^{2}V(\eta_{0})(1+\norm{\ybf^{*}})}{\vep^{2}}}}}.$
 At the $T$-th stage, the algorithm will terminate within $\mcal O\brbra{\log\rbra{{1}/{\vep}}+\rbra{{1}/{\vep}}^{{2}/(1+3\rho)}}$ iterations. Therefore, the total iteration number is of the order $\mcal O\brbra{\rbra{\norm{\ybf^{*}}_{1}+1}\log\brbra{\norm{\ybf^{*}}_{1}+1}\rbra{{1}/{\vep}}^{{2}/(1+3\rho)}}$.

 \blue{For Algorithm~\ref{alg:secant-2}, } 
we note that $T\le T^{\rm{sc}}_\vep$. 
If  $ T_{\ybf^*}\le T^{\rm{sc}}_\vep-1$, following \eqref{eq:upper_bound_sum} and setting $s=s^*$, we have
\begin{equation*}
    \begin{aligned}
    & \sum_{t=1}^{T-1}\min\bcbra{\frac{1}{V(\eta_{t})},\frac{1}{\vep}}^{\frac{2}{1+3\rho}}
    \le  \left(\frac{ |V'(\eta_{0})|}{\vep}\right)^{\frac{2}{1+3\rho}}+T_{\ybf^*}\left(\frac{1}{\vep}\right)^{\frac{2}{1+3\rho}} \\
    & \quad = \Ocal\Big(\left((\norm{\ybf^*}_1+1)\log(\norm{\ybf^*}_1+1)\right)\left(\frac{1}{\vep}\right)^{\frac{2}{1+3\rho}}\Big).
    \end{aligned}
\end{equation*}

If $T_{\ybf^*}> T^{\rm{sc}}_\vep-1$, then setting $s=T^{\rm{sc}}_\vep-1$, we have
$\tsum_{t=1}^{T-1}\min\bcbra{\frac{1}{V(\eta_{t})},\frac{1}{\vep}}^{{2}/(1+3\rho)} \le T^{\rm{sc}}_\vep \brbra{\frac{1}{\vep}}^{{2}/(1+3\rho)}=\Ocal\Big(\log\brbra{\frac{1}{\vep}}\brbra{\frac{1}{\vep}}^{{2}/(1+3\rho)}\Big).$
Furthermore, since $V(\eta_t) \geq \vep /\alpha$, $T^{\rm{sc}}_\vep=\mcal O\Brbra{\min \bcbra{{\frac{\alpha}{-\beta \blue{\bar{V}^{\prime}}}\log\brbra{\frac{\alpha V(\eta_{0})}{-\blue{\bar{V}^{\prime}}\vep}}}, \log\brbra{\frac{\alpha l_{0}}{\varepsilon}\sqrt{\frac{\blue{f^*}-\eta_{0}}{\eta_{1}-\eta_{0}}}}}}$ in Theorem~\ref{thm:convergence-secant-1} and $\blue{\bar{V}^{\prime}}=-\frac{1}{1+\norm{\ybf^*}}$, we have
\begin{equation}
\begin{aligned}\mcal O\Brbra{\sum_{t=1}^{T-1}\log(\frac{1}{V(\eta_{t})})}=\mcal O{\Brbra{\min\bcbra{\frac{\alpha(1+\norm{\ybf^{*}})}{\beta}\log\brbra{\frac{\alpha^{2}V(\eta_{0})(1+\norm{\ybf^{*}})}{\vep^{2}}},\log\brbra{\frac{\alpha^{2}l_{0}}{\varepsilon^{2}}\sqrt{\frac{\blue{f^*}-\eta_{0}}{\eta_{1}-\eta_{0}}}}}}}\end{aligned}
\end{equation}
At the $T$-th stage, the algorithm will end in $\mcal O\Brbra{\log\brbra{{1}/{\vep}}+\brbra{{1}/{\vep}}^{{2}/(1+3\rho)}}$ iterations. Therefore, the total iteration number is of the order $\mcal O\Brbra{\min \bcbra{\rbra{\norm{\ybf^*}_1+1}\log\rbra{\norm{\ybf^*}_1+1}, \log\brbra{{1}/{\vep}}} \cdot \brbra{{1}/{\vep}}^{{2}/(1+3\rho)}}$.

\section{The Initialization Phase}\label{subsec:initialization}

We develop routines to find the initial point. We first give some intuitions on how to find the initial level $\eta_0$ for the fixed point iteration, which, together with an initial point $\tilde{\xbf}^0$, shall satisfy the condition:
\begin{equation}\label{eq:strict-infeasible}
\eta^0=f(\tilde{\xbf}^0) < \blue{f^*}. 
\end{equation}
Let $\bar{f}:=f(\bar{\xbf})$.
We immediately observe that $f(\bar{\xbf})\le \blue{f^*}$.
If $\bar{\xbf}$ is feasible to the set constraint $\{g(\xbf)\le \zero\}$, then we know $\bar{\xbf}$ is optimal to the original problem~\eqref{pb:func-constraint}, and $f(\bar{\xbf})=f^*$. Otherwise, $\bar{\xbf}$ must satisfy \eqref{eq:strict-infeasible}.

In practice, finding an exact minimizer $\bar{\xbf}$ can be highly challenging. To address this, we develop an initialization scheme in Algorithm~\ref{alg:phase-I}. 
First, we compute an $\varepsilon$-accurate minimizer $\tilde{\xbf}^0\in\Xcal$: $f(\tilde{\xbf}^0)-\bar{f}\le \vep$.
If $\max_{1\le i\le m} \{g_i(\tilde{\xbf}^0)\}\le \vep$, then $\tilde{\xbf}^0$ is already an $\vep$-optimal solution to problem~\eqref{pb:func-constraint}.
Otherwise, we approximately solve the fixed-level problem $\min_{\xbf\in\Xcal}\, v(\xbf, f(\tilde{\xbf}^0)).$
The following theorem guarantees that by solving the fixed-level problem, Algorithm~\ref{alg:phase-I} will either produce a near-optimal solution and terminate or it will generate a valid solution to initiate the main root-finding phase.
\begin{theorem}
When Algorithm~\ref{alg:phase-I} enters line~\ref{step:call_apl}, 
the number of gap reduction iterations is bounded by 
    \begin{equation}\label{eq:terminate_iteration}
            S_0 + \frac{1}{1-\gamma^{2/(1+3\rho)}}\Big(\frac{2^{\rho+2}3^{(1-\rho)/2}\bar{M}D_\Xcal^{(\rho+1)/2}\cdot \alpha}{(1+\rho)(\alpha-1)}\min\{\frac{1}{V(\eta_0)},\frac{1}{\vep}\}\Big)^{\frac{2}{1+3\rho}},
    \end{equation}
where  $S_0= \max\left\{0, \left\lceil\log_{1/\gamma} \Brbra{\frac{\alpha(\bar{u}_0-\bar{l}_0)}{(\alpha-1)}\min\{\frac{1}{V(\eta_0)},\frac{2(\alpha-1)}{\alpha \vep}\}}\right\rceil\right\}.$
   The solution returned by Algorithm~\ref{alg:phase-I} will either be an  $\vep$-optimal solution or satisfy the condition in \eqref{eq:strict-infeasible}. 
\end{theorem}
\begin{proof}
Algorithm~\ref{alg:phase-I} enters line~\ref{step:call_apl} when $\tilde{\xbf}^0$ fails to be $\vep$-optimal to \eqref{pb:func-constraint}. 
By setting $\theta = \frac{1}{2}$ in Theorem~\ref{thm:complexity-APL}, we conclude that the algorithm will terminate in at most the number of iterations specified in~\eqref{eq:terminate_iteration}.
APL will terminate in two conditions. 
If $\tilde{u}\le \vep$, then we have $
    f(\tilde{\xbf}^*)-\eta_0\le \vep, \ \max_{1\le i\le m} g_i(\tilde{\xbf}^*)\le \vep.$ 
    Since $\bar{f}\leq f^*$,  we have $\eta_0-\blue{f^*} = f(\tilde{\xbf}^0) -f^* \le f(\tilde{\xbf}^0) -\bar{f}\le \vep.$
    It follows that $\tilde{\xbf}^*$ is a $2\vep$-optimal solution of \eqref{pb:func-constraint}. 
    If $\tilde{u}> \vep$, we must have $\tilde{u}-\tilde{l}\le (\alpha-1)/\alpha \tilde{u}$ and hence $\tilde{l}>0$. It follows that $V(\eta_0)\ge\tilde{l}>0=V(\blue{f^*})$. 
    Since $V(\cdot)$ is monotonically decreasing (from Proposition~\ref{prop:value-function}),  we have $\eta_0<\blue{f^*}$. We conclude that the assumptions of the root-finding are satisfied.  

\end{proof}

\begin{algorithm}[h!]
    \caption{The Initialization Phase}\label{alg:phase-I}
    \KwIn{$\alpha$, $\vep$;}
    Compute an $\vep$-optimal solution $\tilde{\xbf}^0\in\Xcal$: $f(\tilde{\xbf}^0)-\bar{f}\le \vep$ and set $\eta_{0}=f({\tilde\xbf}^{0})$\;
    \If{$\max_i\{g_i(\tilde{\xbf}^0)\}\le \vep$}{
    \KwRet{$(\tilde{\xbf}^0, \eta_{0}, -\infty, \True)$;}   \tcp{A near-optimal point found} 
    }
    Set $l_0=\min_{\xbf\in\Xcal} v_\ell (\xbf, \tilde{\xbf}^0, \eta_0)$\label{step:initial-next}\; 
    Compute $(\tilde{\xbf}^*, \tilde{l})=\Acal(\tilde{\xbf}^0, l_0, \eta_0, \frac{1}{2}, \alpha)$  and set $\tilde{u}=v(\tilde{x}^*, \eta_0)$\label{step:call_apl}\;
    \eIf{$\tilde{u} \le \vep $}{
        \KwRet{$(\tilde{\xbf}^*, \eta_{0}, -\infty, \True)$;}  \tcp{A near-optimal point found} 
    }{
    \KwRet{$(\tilde{\xbf}^*, \eta_{0}, \tilde{l},  \False)$;} \tcp{Root-finding required}
    }
\end{algorithm}
\begin{remark}
Our analysis does not include the complexity of computing $\tilde{\xbf}^0$. 
Since minimizing $f(\xbf)$ over $\Xcal$ does not involve function constraints, 
it is much easier to solve than problem~\eqref{pb:func-constraint}. 
Specifically, when applying the universally optimal accelerated method~\citep{nesterov2015universal} or the original APL method~\citep{lan2015bundle}, one can obtain $\tilde{\xbf}^0$ in $\Ocal(\vep^{-2/(1+3\rho)})$ iterations. Therefore, this step does not affect the complexity order of Algorithm~\ref{alg:phase-I}. 
\end{remark}

   
\end{document}